\documentclass[11pt,a4paper]{amsart}
\usepackage{amsmath,amsthm,amssymb}
\usepackage{comment}
\newcommand{\mychoice}[3]{#1
}
\mychoice{
\newcommand{\plabel}[1]{ \label{#1}}
\newcommand{\gbibitem}[1]{ \bibitem{#1}}
\newcommand{\snewpage}{}

\specialcomment{commentx}{}{}\excludecomment{commentx}
\specialcomment{commenty}{}{}\excludecomment{commenty}
\specialcomment{commentz}{}{}

}{
\usepackage{xcolor}
\newcommand{\plabel}[1]{ \label{#1}\rlap{\smash{${}^{^{[#1]}}$}}}
\newcommand{\gbibitem}[1]{ \bibitem{#1}\rlap{\smash{${}^{^{[#1]}}$}}}
\newcommand{\snewpage}{\newpage}

\newenvironment{commentx}{\color{magenta} }{\color{black} }
\newenvironment{commenty}{\color{blue} }{\color{black} }

}{
\usepackage{xcolor}
\newcommand{\plabel}[1]{ \label{#1}}
\newcommand{\gbibitem}[1]{ \bibitem{#1}}
\newcommand{\snewpage}{}

\specialcomment{commentx}{}{}\excludecomment{commentx}

\specialcomment{commentz}{}{}\excludecomment{commentz}

}
\DeclareMathOperator{\im}{im}
\DeclareMathOperator*{\subs}{subs}
\newcommand{\ext}{\mathrm{ext}}

\DeclareMathOperator{\nonass}{n-a}

\DeclareMathOperator{\Lie}{Lie}
\DeclareMathOperator{\Der}{Der}
\DeclareMathOperator{\proj}{pr}
\DeclareMathOperator{\Words}{Words}
\DeclareMathOperator{\ad}{ad}
\DeclareMathOperator{\id}{id}
\DeclareMathOperator{\cmr}{cmr}
\newcommand{\bo}{\boldsymbol}

\newcommand{\botimes}{{\textstyle\bigotimes}}
\newcommand{\bodot}{{\textstyle\bigodot}}
\newcommand{\leaveout}[1]{}

\newcommand{\ass}{\mathrm{assoc}}
\newcommand{\lw}{\mathrm{Lw}}
\theoremstyle{definition}
\newtheorem{point}{}[section]
\newtheorem{remark}[point]{Remark}

\theoremstyle{plain}
\newtheorem{prop}[point]{Proposition}
\newtheorem{lemma}[point]{Lemma}
\newtheorem{theorem}[point]{Theorem}
\newtheorem{cor}[point]{Corollary}
\newcommand{\qedremark}{  \renewcommand{\qedsymbol}{$\triangle$} \qed \renewcommand{\qedsymbol}{$\Box$}}

\newcommand{\eqed}{
\pushQED{\qed}
\qedhere
\popQED
}
\newcommand{\marginextend}[1]{ \addtolength{\oddsidemargin}{-#1}  \addtolength{\evensidemargin}{-#1}\addtolength{\textwidth}{#1}\addtolength{\textwidth}{#1}}
\newcommand{\updownextend}[1]{ \addtolength{\topmargin}{-#1}  \addtolength{\textheight}{#1}
\addtolength{\textheight}{#1}}
\marginextend{1.25cm}
\updownextend{0cm}
\allowdisplaybreaks[3]
\begin{document}
\title{Some proofs of the Poincar\'e--Birkhoff--Witt theorem and related matters}
\date{\today}
\author{Gyula Lakos}
\email{lakos@renyi.hu}
\address{Alfréd Rényi Institute of Mathematics}
\begin{abstract}
This expository paper focuses on free Lie $K$-algebras and the basic PBW theorem.
We argue in various ways that the basic PBW theorem is a quite close consequence of the Magnus--Witt theorems concerning free Lie algebras.
\end{abstract}
\maketitle

\section*{Introduction}\plabel{sec:intro}

The objective of this paper is to consider alternative proofs of the basic Poincar\'e--Birkhoff--Witt theorem
 using universal algebraic arguments.
In this present version of the paper, under `related matters' we will chiefly mean free Lie algebras
 and an array of techniques applicable to them.
First of all, we should make clear that the basic PBW theorem is not a particularly difficult theorem,
 and investing some time in reading a proof of it as in, say, Jacobson \cite{J} can deal with the matter well.
Nevertheless, most proofs of the PBW theorem are either computational or combinatorial.
Not in an extreme manner, but either the computations or the combinatorics can leave some unease that
 we were just ``too lucky'' to have a proof.
On the other hand, there are people who find universal algebraic arguments relatively effortless.
This paper is written for them.
We will argue that the basic PBW theorem in a natural development of the basic Magnus--Witt theorems
of Magnus \cite{MM} and Witt \cite{W} concerning free Lie algebras,
 and, in particular, one obtains the PBW theorem in this venue quite naturally.
This is certainly well-known to the experts, and
 they make this point in various ways, see \v{S}ir\v{s}ov \cite{S2}, Lothaire \cite{Lot}, Reutenauer \cite{R}.
Here multiple viewpoints are collected, except with certain changes.

Let us first review the statement of the PBW theorem, similarly to \cite{LLX} but in slightly more detail.

\textbf{The local formulation of the PBW theorem.}
Assume that $K$ is a unital commutative ring, and $\mathfrak g$ is a $K$-module with a compatible Lie-ring structure;
 i.~e.~$\mathfrak g$ is a Lie $K$-algebra (also called: Lie ring over $K$).
The universal enveloping algebra $\mathcal U\mathfrak g$ is the free $K$-algebra
 $\mathrm F_K[\mathfrak g]\simeq \bigotimes\mathfrak g\equiv \bigoplus_{n=0}^\infty\botimes^n\mathfrak g$
 factorized by the ideal $J\mathfrak g$ generated by the
 elements $X\otimes Y-Y\otimes X-\boldsymbol[X,Y\boldsymbol]$, the tensor products are taken over $K$.
Let $\boldsymbol m:\bigotimes\mathfrak g\rightarrow\mathcal U\mathfrak g$ denote this canonical homomorphism.
The enveloping algebra is naturally filtered by $\mathcal U^n\mathfrak g=\bo m\left(\botimes^{\leq n}\mathfrak g\right)$, and the
 construction implies the existence of natural surjective maps
 $\bo m^{(n)}: \bodot^n\mathfrak g\rightarrow \mathcal U^n\mathfrak g/\mathcal U^{n-1}\mathfrak g$.
The (local form of the) Poincar\'e--Birkhoff--Witt theorem, whenever it holds, states  that the maps $\bo m^{(n)}$ are isomorphisms.
This theorem is known to hold in the following cases:

(i) $K$ is a field, or, more generally,  $\mathfrak g$ is a free $K$-module, or,
 more generally, $\mathfrak g$ is a direct sum of cyclic $K$-modules
 (Poincar\'e \cite{P}: $K$ is a field, $\mathbb Q\subset K$, cf.~Ton-That, Tran \cite{TT};
 Birkhoff \cite{B}, Witt \cite{W}: $K$ is a field, but their methods work more generally);

(i') $K$ is a principal ideal domain (Lazard \cite{L}) or just a Dedekind domain (Cartier \cite{C});
 also see Higgins \cite{H} for further results in this direction;

(ii) $\mathbb Q\subset K$ (Cohn \cite{CC});

(ii') $\frac12\in K$ but $[\mathfrak g,[\mathfrak g,\mathfrak g]]=0$ (Nouaz\'e, Revoy \cite{NR}),
or  $\frac1{k!}\in K$ but $(\ad\mathfrak g)^{k}\mathfrak g=0$ (\cite{LLX}) ;
\\
cf.~Grivel \cite{G} for a review. But there are counterexamples
 (see \v{S}ir\v{s}ov \cite{S1}, Cartier \cite{C}, Cohn \cite{CC}).
The most general approach is of Higgins \cite{H}, cf.~Revoy \cite{Re}.


\textbf{Global formulations of the PBW theorem.}
In practice, mostly cases  (i) and (ii) are considered but in global form:

(i) If $\mathfrak g$ is a sum of cyclic $K$-modules, then we can choose a basis
 $\{g_\alpha\,:\,\alpha\in A\}$, and an ordering $\leq$ of $A$.
Then let $\bigotimes_\leq\mathfrak g$ be the submodule of $\bigotimes\mathfrak g$ spanned by
$g_{\alpha_1}\otimes \ldots\otimes g_{\alpha_n}$ with $\alpha_1\leq\ldots\leq\alpha_n$.
Then the ``basic'' version of the PBW theorem states that
\begin{equation*}\boldsymbol m_\leq:\textstyle{\bigotimes_\leq\mathfrak g}\rightarrow\mathcal U\mathfrak g\end{equation*}
 (a restriction of $\boldsymbol m$) is an isomorphism of $K$-modules.
(In fact, what we really need, in general, is a choice function $\mathfrak c$, which transforms any
 finite multiset $\bo\alpha$ from $A$ into an ordered word $\alpha_1 \ldots \alpha_n$.
Then the statement is that is that the corresponding map
 $\boldsymbol m_{\mathfrak c}:\textstyle{\bigotimes_{\mathfrak c}\mathfrak g}\rightarrow\mathcal U\mathfrak g$
 is an isomorphism.
For the sake of simplicity, we will use the ordered version.)
As the natural map $\bigotimes_\leq^n\mathfrak g\rightarrow\bodot^n\mathfrak g$ is an isomorphism,
 it is easy to see that this $\boldsymbol m_\leq$ is necessarily surjective,
 and its injectivity is equivalent to the local version of PBW theorem; one can see these inductively on the filtration.
The surjectivity can be interpreted in a semi-algorithmical way (independently from the local formulation):
Whenever we have an expression in $\mathcal U\mathfrak g$ (as an image of $\bo m$),
 we can rearrange the formally top nonarranged degree term into an image of $\bo m_\leq$
 at the cost of generating formally lower order terms.
Then we repeat this in formally lower orders. We call this as the ``basic rearrangement procedure''.

(ii) If $\mathbb Q\subset K$, then one can consider the submodule $\bigotimes_\Sigma\mathfrak g$ of $\bigotimes\mathfrak g$.
This  submodule can be interpreted either as the submodule of elements invariant under permutations
 in the order of tensor product or as the span of the elements
 $a_1\otimes_\Sigma \ldots\otimes_\Sigma a_n=\frac1{n!}\sum_{\sigma\in\Sigma_n}g_{\sigma(1)}\otimes\ldots\otimes g_{\sigma(n)}$.
Then the ``symmetric'' version of the PBW theorem states that
\begin{equation*}\boldsymbol m_\Sigma:\textstyle{\bigotimes_\Sigma\mathfrak g}\rightarrow\mathcal U\mathfrak g\end{equation*}
 (a restriction of $\boldsymbol m$) is an isomorphism of $K$-modules.
As the natural map $\bigotimes_\Sigma^n\mathfrak g\rightarrow\bodot^n\mathfrak g$ is an isomorphism,
 it is easy to see that this $\boldsymbol m_\leq$ is necessarily surjective,
 and its injectivity is equivalent to the local version of PBW theorem;
 one can see these inductively on the filtration.
The surjectivity can be interpreted in a semi-algorithmical way (independently from the local formulation):
We can demonstrate the surjectivity of $\boldsymbol m_\Sigma$
 using the ``symmetric rearrangement procedure'' in $\mathcal U\mathfrak g$,
 i.~e.~symmetrizing in the formally top nonarranged degree term at the cost of generating lower order terms,
 and repeating the process in formally lower orders.

In practice, one typically starts with the global cases (i), (ii), and then proceeds further to the general local ones.
The Witt--Lazard version of the proof of the PBW theorem can be adapted to cases (i) and (ii) alike.
Such a formulation is given in Appendix \ref{sec:lazard}.
In terms of the basic PBW theorem, i.~e.~case (i), Jacobson \cite{J} writes up the somewhat terse account of Birkhoff \cite{B};
 but the result is also not much different from the ideas of Witt \cite{W}.
Some other classical proofs are also commented in Appendix \ref{sec:lazard}.

Nevertheless, it can be said  that case (ii) is easier than case (i):
The proof of  (ii) is just one step away from the statement that the PBW theorem holds
 for free Lie algebras $\mathrm F_{\mathbb Q}^{\Lie}[X_\lambda : \lambda\in\Lambda]$;
 while  (i) sort of needs $\mathrm F_{\mathbb Z}^{\Lie}[X_\lambda : \lambda\in\Lambda]$ for this.
By this, we arrive to

\textbf{The free case.}
One of the few cases where $\mathcal U\mathfrak g$ is easy to describe is
 when $\mathfrak g$ is the free Lie $K$-algebra $\mathrm F_K^{\Lie}[X_\lambda : \lambda\in\Lambda]$.
Then $\mathcal U\mathrm F_K^{\Lie}[X_\lambda : \lambda\in\Lambda]$ is naturally isomorphic to the noncommutative
 polynomial algebra $\mathrm F_K[X_\lambda : \lambda\in\Lambda]$.
This holds for purely universal algebraic reasons. Indeed, $\mathrm F_K^{\Lie}[X_\lambda : \lambda\in\Lambda]$ evaluates
 by commutators, which defines a map
 $\iota:\mathrm F_K^{\Lie}[X_\lambda : \lambda\in\Lambda]\rightarrow \mathrm F_K[X_\lambda : \lambda\in\Lambda]$,
 which, by the universality of the enveloping algebra, gives rise to a map
 $\mathcal U \iota:\mathcal U\mathrm F^{\Lie}[X_\lambda : \lambda\in\Lambda]\rightarrow \mathrm F_K[X_\lambda : \lambda\in\Lambda]$.
Here the class $[X_\lambda]_{\mathcal U}$ maps to $X_\lambda$.
Comparing this to the universality of $\mathrm F_K[X_\lambda : \lambda\in\Lambda]$, we see that $\mathcal U \iota$ is an isomorphism.

What is not obvious is that
 $\iota:\mathrm F_K^{\Lie}[X_\lambda : \lambda\in\Lambda]\rightarrow
 \mathrm F_K[X_\lambda : \lambda\in\Lambda]\simeq\mathcal U\mathrm F_K^{\Lie}[X_\lambda : \lambda\in\Lambda]$
 is an inclusion.
This is the Magnus--Witt theorem of on the representability of free Lie $K$-algebras
 (cf.~Magnus \cite{MM} and Witt \cite{W}; the critical case is $K=\mathbb Z$).
But this is a just the weak version of the PBW theorem ($\bo m^{(1)}$ level).
Now, the PBW theorem does, indeed, hold in the case when $\mathfrak g$ is a free Lie $K$-algebra;
 because it is true that free Lie $K$-algebras are free $K$-modules.
This latter fact is not trivial if $K$ is not a field, but it is the Magnus--Witt' theorem on freely generatedness,
 which, in turn, can be derived from the  Magnus--Witt theorem for $K=\mathbb Z$.
This circle can be resolved by applying the basic theorem PBW  for $K=\mathbb Z$,
 using that we have a multigrade-wise finitely generated $\mathbb Z$-module,
 which must be a direct sum of cyclic $\mathbb Z$-modules.

While thus the free case, ultimately, fits to the global case (i), free Lie $K$-algebras have a different, more combinatorial flavour.
Firstly, the Magnus--Witt theorems can be proven quite easily, by
 elimination; by the Lie algebraic Lazard--Shirshov method of \v{S}ir\v{s}ov \cite{S0} and Lazard \cite{L2}
 (cf.~Bourbaki \cite{BX}),
 or by the original noncommutative polynomial elimination of Magnus \cite{MM}
 (cf.~Magnus, Karrass, Solitar \cite{MKS}, Labute \cite{Lab}).
Secondly, quite explicit bases for $\mathrm F_K^{\Lie}[X_\lambda : \lambda\in\Lambda]$ can be written down.
See M. Hall \cite{Ha} for the classical Hall bases (following ideas of P. Hall \cite{Ha0}  and Magnus \cite{MM});
 Chen, Fox, Lyndon \cite{CFL} and \v{S}ir\v{s}ov \cite{S2} for the Lyndon--Shirshov basis;
 \v{S}ir\v{s}ov \cite{S3} and Michel \cite{Mch} and Viennot \cite{V} for their common generalization, the Hall--Shirshov bases;
 these are the most well-known ones.
Utilizing this machinery, the Magnus--Witt theorems also follow immediately.
Thirdly, the above mentioned bases are also well-adapted to the proof of PBW theorem for free Lie algebras due to
 some natural factorization properties.
This was shown by Sch\"utzenberger \cite{Sc}, \v{S}ir\v{s}ov \cite{S2}, Viennot \cite{V} with respect to bases above
(cf. also Lothaire \cite{Lot}).
A detailed discussion of these  can also be found in Reutenauer \cite{R};
 a shorter account directed to general case of Hall--Shirshov bases is given by Reutenauer \cite{RR}.
Fourthly, and finally, the PBW theorem  for free Lie algebras implies the  PBW theorem cases (i) and (ii).
For the basic case (ii), this is already observed by \v{S}ir\v{s}ov \cite{S2}; and a simple and modern
 presentation with the emphasis of the combinatorial use Lyndon--Shirshov words is given by Lothaire \cite{Lot}.

We will approach the PBW theorem in the spirit of the latter paragraph, with an emphasized use of free Lie algebras.
The style itself is not completely introductory but with sufficiently mathematical thinking
 one can read it for his pleasure, and can decide that the arguments given are sufficiently elementary or not.
For a broader view on the PBW theorem, and aspects not considered here we refer to
 Grivel \cite{G},  Shepler, Witherspoon \cite{SW}, Bremner, Dotsenko \cite{BD}.

\snewpage

\textbf{Outline of content.}
In Sections \ref{sec:free1} and \ref{sec:free2}, we include some elementary facts regarding free Lie $K$-algebras,
 and we show how to prove the Magnus--Witt theorems  using a sufficiently strong version of the PBW theorem itself
 or by Lie algebraic elimination or noncommutative polynomial elimination.
(This follows  Magnus \cite{MM}, Witt \cite{W} and \v{S}ir\v{s}ov \cite{S0}, Lazard \cite{L2}.)
Not everything explained here must be used in what follows, but it is useful to obtain a more complete picture.
In Section \ref{sec:via}, we discuss the PBW theorem for free Lie algebras over $\mathbb Q\subset K$,
 and the general free case is also discussed modulo the Magnus--Witt theorem.
(Symmetrization is already a theme in Poincaré \cite{P}.)
In Section \ref{sec:freeMW}, we explain how to obtain PBW theorem for free Lie algebras
 from the Magnus--Witt theorems and the coproduct structure.
(The coproduct is used by  Milnor, Moore \cite{MiMo}, but, beside the use of the Magnus--Witt theorems,
 a  difference is that we use the  collective symmetric coproduct in order to deal with the case $\mathbb Q\not\subset K$  directly.)
In Section \ref{sec:freedirect}, a version of Shirshov's argument for free Lie algebras  is produced.
In Section \ref{sec:freefrom}, it is explained how the free case implies the basic PBW theorem.
(These last two sections are   in the spirit of \v{S}ir\v{s}ov \cite{S2} and Lothaire \cite{Lot}.)

\textbf{Acknowledgements.}
 The author thanks Bal\'azs Csik\'os and Márton Naszódi.

\snewpage
\section{About free Lie algebras. Version 1 (Abstract version)}\plabel{sec:free1}
\textbf{Gradedness and its stability.}
Free Lie $K$-algebras (or any other kinds of free algebras) do not really require specific constructions.
Nevertheless, it  is very useful to have various theorems which provide some control over them.
A simple but extremely useful property is that free Lie $K$-algebras have natural gradings.
This statement is obvious in case of the free non-associative algebras (magmas, brace algebras) $\mathrm F^{\nonass}_K$
 and free  associative algebras $\mathrm F_K$, but it requires some care in free Lie algebras.
Let us think about the free Lie $K$-algebra $\mathrm F_K^{\Lie}[X_\lambda :\lambda\in\Lambda]$
 as the free nonassociative $K$-algebra $\mathrm F^{\nonass}_K[X_\lambda :\lambda\in\Lambda]$
 factorized further by the $K$-submodule (ideal) $\mathrm I_K^{\Lie}[X_\lambda :\lambda\in\Lambda]$.
Additively, $\mathrm F^{\nonass}_K[X_\lambda :\lambda\in\Lambda]$ is just the
 free $K$-module generated by the $\bo[,\bo]$-monomials of the $X_\lambda$.
\begin{lemma}\plabel{lem:freelie}
$\mathrm I^{\Lie}_K[X_\lambda :\lambda\in\Lambda]$ is generated by the elements
\begin{enumerate}
\item[(F1)] $k M( \bo[Z_1,Z_1\bo],X_{\lambda_1},\ldots, X_{\lambda_s})$ ,
\item[(F2)] $k M( \bo[Z_1,Z_2\bo],X_{\lambda_1},\ldots, X_{\lambda_s})+kM( \bo[Z_2,Z_1\bo],X_{\lambda_1},\ldots, X_{\lambda_s}) $ ,
\item[(F3)] $k M( \bo[\bo[Z_1,Z_2\bo],Z_3\bo],X_{\lambda_1},\ldots, X_{\lambda_s})$ $+kM( \bo[\bo[Z_2,Z_3\bo],Z_1\bo],X_{\lambda_1},\ldots, X_{\lambda_s})$

$+kM( \bo[\bo[Z_3,Z_1\bo],Z_2\bo],X_{\lambda_1},\ldots, X_{\lambda_s})
 $;\end{enumerate}
 where $Z_1,Z_2,Z_3$ are monomials of the $X_\lambda$, and $M(\ldots)$ is a $\bo[,\bo]$-bracketing
 with $s+1$ many positions (but not necessarily in the indicated order), and $k\in K$.
\begin{proof}
Such elements are clearly in the ideal $\mathrm I^{\Lie}_K$.
Conversely, whenever we take elements from  $\mathrm F^{\nonass}_K$ and
 apply the Lie-identities, then they expand to sums of cases (F1)--(F3) with trivial $M$.
(Note that case (F2) cannot be omitted.)
Thus the primary relations (coming form the Lie-identities) are generated. The secondary relations
 (coming from $x\sim y\rightarrow \bo[x,z\bo]\sim \bo[y,z\bo], \bo[z,x\bo]\sim \bo[z,y\bo]$ are
 also generated due to linearity and that nontrivial $M$ are allowed.
\end{proof}
\end{lemma}
\begin{prop}\plabel{prop:freelie}
(a) $\mathrm F^{\Lie}_K[X_\lambda :\lambda\in\Lambda]$ is multigraded by the number of various variables.

(b) The structure in a given multigrade depends only on its multiplicity structure
 (independently of the presence of other variables, etc.).
\begin{proof}
(a) The multigradedness will be inherited from $\mathrm F^{\nonass}_K[X_\lambda :\lambda\in\Lambda]$, because
 the relations from Lemma \ref{lem:freelie} are multigrade-homogeneous.

(b) The structure of $\mathrm I^{\Lie}_K[X_\lambda :\lambda\in\Lambda]$ also depends only
 on the multiplicity pattern.
\end{proof}
\end{prop}
\begin{cor}\plabel{cor:freelie}
In any  multigrade, corresponding to finite multiset of the $X_\lambda$,
 $\mathrm F^{\Lie}_K[X_\lambda :\lambda\in\Lambda]$ is generated by finitely many monomials.
\begin{proof}
Every finite multiset of the $X_\lambda$ can be bracketed only in finitely many ways, so finitely
 generatedness is true even in $\mathrm F^{\nonass}_K[X_\lambda :\lambda\in\Lambda]$.
\end{proof}
\end{cor}
The following uses not even the grading, but it is compatible with it:
\begin{prop}\plabel{prop:coeff}
$\mathrm F^{\Lie}_K[X_\lambda:\lambda\in\Lambda]\simeq \mathrm F^{\Lie}_{\mathbb Z}[X_\lambda:\lambda\in\Lambda]\otimes K$ naturally.
\begin{proof}
The natural map
 $\mathrm F^{\Lie}_{\mathbb Z}[X_\lambda:\lambda\in\Lambda]\otimes K\rightarrow\mathrm F^{\Lie}_K[X_\lambda:\lambda\in\Lambda]$
 yields a surjective morphism of Lie $K$-algebras, which can be inverted back due to the universal property of the range.
\end{proof}
\end{prop}

\snewpage
\textbf{The Magnus--Witt theorems.} These are the following statements:
\begin{theorem}\plabel{prop:WMag}  (The Magnus--Witt theorem.)
$\mathrm F^{\Lie}_K[X_\lambda:\lambda\in\Lambda]$ embeds to the noncommutative polynomial algebra
 $\mathrm F_K[X_\lambda:\lambda\in\Lambda]$ by the commutator-evaluation.
\end{theorem}
Note that Theorem \ref{prop:WMag} implies Proposition \ref{prop:freelie} (and thus Corollary \ref{cor:freelie}).
\begin{theorem}\plabel{prop:CMag} (The Magnus--Witt' theorem.)
$\mathrm F^{\Lie}_K[X_\lambda:\lambda\in\Lambda]$ is a free $K$-module (multigrade-wise).
[Strong version: With a set of $\bo[,\bo]$-monomials  without reference to $K$ as a basis.]
\end{theorem}
The following terminology is not historical but it sort of fits to the pattern:
\begin{theorem}\plabel{prop:XCMag} (The Magnus--Witt'' theorem.)
The image $\mathrm F^{\cmr}_K[X_\lambda:\lambda\in\Lambda]$ of $\mathrm F^{\Lie}_K[X_\lambda:\lambda\in\Lambda]$
 by the commutator-evaluation can be directly complemented in $\mathrm F_K[X_\lambda:\lambda\in\Lambda]$ (multigrade-wise).
[Strong version: By a free $K$-module, with
 a set of products of $[,]$-monomials without reference to $K$ as a basis.]
\end{theorem}
In what follows, the `Magnus--Witt (MW) theorem' refers to Theorem \ref{prop:WMag}, but the
 `Magnus--Witt (MW) theorems' refers to Theorems \ref{prop:WMag}, \ref{prop:CMag}, \ref{prop:XCMag} collectively.
By the Magnus--Witt' theorem and Magnus--Witt'' theorem we mean the ``weak'' versions of the statements
 unless said otherwise.
(But, as we will see, not even the strong versions are particularly strong.)

\textbf{The Magnus--Witt theorems via the PBW theorem.}
The gradedness allows to apply the PBW theorem (for sum of cyclic submodules) to obtain the (weak versions of) Magnus--Witt theorems.
\begin{proof}[Proof of the MW' theorem via the PBW theorem.]
Assume $K=\mathbb Z$.
Then $\mathrm F^{\Lie}_{\mathbb Z}[X_1,\ldots,X_n]$ is a finitely generated $\mathbb Z$-module in every multigrade,
 thus it is a sum of cyclic $\mathbb Z$-modules.
Then the PBW theorem  (for sums of cyclic submodules) can be applied to show that
 $\mathrm F^{\Lie}_{\mathbb Z}[X_1,\ldots,X_n]$ embeds into $\mathcal U\mathrm F^{\Lie}_{\mathbb Z}[X_1,\ldots,X_n]\simeq \mathrm F^{}_{\mathbb Z}[X_1,\ldots,X_n]$.
It is immediate that restricted to $\mathrm F^{\Lie}_{\mathbb Z}[X_1,\ldots,X_n]$ it is given by the  commutator evaluation.
As the image of $\mathrm F^{\Lie}_{\mathbb Z}[X_1,\ldots,X_n]$ has no torsion, so
$\mathrm F^{\Lie}_{\mathbb Z}[X_1,\ldots,X_n]$ is a free $\mathbb Z$-module in every multigrade.
By Proposition \ref{prop:coeff}, the MW' theorem for general $K$ holds.
(Or directly: we observe that additively $\mathrm F^{\nonass}_{\mathbb Z}[X_1,\ldots,X_n]
 \simeq (\text{a free $\mathbb Z$-module})\oplus\mathrm I^{\Lie}_{\mathbb Z}[X_1,\ldots,X_n]$
(in every multigrade).
This decomposition structure survives by tensoring with $K$, so general case follows.)
\end{proof}
\begin{proof}[Proof of the MW and MW'' theorems via the PBW theorem.]
The MW and MW'' the\-o\-rems  follow immediately by considering the PBW theorem (for free modules).
\end{proof}

\textbf{General observations regarding the Magnus--Witt theorems.}
The case $K=\mathbb Z$ in the Magnus-Witt theorems implies the general case by relatively simple arguments.
Indeed, in the cases of MW' and MW'' the implication is immediate.
Regarding the MW theorem, it is sufficient to see that the natural map
 $\mathrm F^{\cmr}_{\mathbb Z}[X_\lambda:\lambda\in\Lambda]\otimes K\rightarrow\mathrm F^{\cmr}_K[X_\lambda:\lambda\in\Lambda]$
 is an isomorphism, but this follows from the MW'' theorem.
It can be argued that only \v{S}ir\v{s}ov \cite{S2} addresses the case of general $K$ as such,
 but the case $K=\mathbb Z$ is done by Magnus \cite{MM} and Witt \cite{W}.
The MW' and MW'' theorems are trivial if $K$ is a field.
The MW theorem follows from (an appropriate formulation of) the Dynkin--Specht--Wever lemma if $\mathbb Q\subset K$.
Moroever, if $K$ is a field, then a less sophisticated version of the PBW theorem can be used.
In general, however, the technique of elimination provides a simple proof for
 the MW theorem (and also to MW' and MW'' theorems):

\snewpage

\textbf{Elimination by Lie algebraic derivations.}
It is easy to see that derivations $D$ of $\mathrm F^{\nonass}_K[X_\lambda :\lambda\in\Lambda]$ are determined by arbitrary prescriptions
$D(X_\lambda)=P^{\nonass}_\lambda$.
Then Lemma \ref{lem:freelie} implies easily that these derivations descend to  $\mathrm F^{\Lie}_K[X_\lambda :\lambda\in\Lambda]$.
In particular, its derivations are also given by  arbitrary prescriptions $D(X_\lambda)=P_\lambda^{\Lie}$.
We can apply this in the setting when $\Lambda$ is a the set of words $\Words(L)$ on an alphabet $L$,
$\mathfrak g=\mathrm F^{\Lie}_K[X_\lambda :\lambda\in\Lambda]$,
 and the derivations
$D_\nu $ $(\nu\in \Words(L))$ are given by $D_\nu(X_\lambda)=X_{\nu\lambda}$.
It is easy  to see that these derivations \textit{linearly} generate the left word derivation
Lie $K$-algebra $\Der^{\lw}_K \mathfrak g$.
(The commutator is $[D_\nu,D_\xi]=D_{\nu\xi}-D_{\xi\nu}$.)
Then we can define  the extension $(\Der^{\lw}_K \mathfrak g)\ltimes\mathfrak g$ such that $\bo[D,x\bo]=D(x)$.

Elimination theory uses the following more specific setup:
Assume $\Xi=A\,\dot\cup\,B$, where $A,B\neq\emptyset$.
Let $\mathfrak h=\mathrm F^{\Lie}_K[X_\xi :\xi\in\Xi]$.
Then $\mathfrak h=\mathfrak h_{\deg_B=0}\ltimes\mathfrak h_{\deg_B>0}$
 (inner semidirect sum), where $\mathfrak h_{\deg_B=0}$ naturally identifies with
 $\mathrm F^{\Lie}_K[X_\alpha :\alpha\in A]$, and $\mathfrak h_{\deg_B>0}$ is the part with multidegrees when
 some variable $X_\beta$ with $\beta\in B$ has multiplicity at least $1$.
As $\mathfrak h$ is a free Lie algebra, there is a homomorphism $\theta_{A,B}$ to the extended Lie algebra of the
 previous paragraph with $L=A\cup B$ given by $\alpha\mapsto D_\alpha$  $(\alpha\in A)$
 and  $\beta\mapsto X_\beta$  $(\beta\in B)$.
Then, using standard Lie rules (absorbing the $D_\alpha$), it is easy to see that $\theta_{A,B}(\mathfrak h_{\deg_B>0})$
 is generated by $X_{\alpha_1\ldots\alpha_n\beta}$.
Hence the image is the free Lie algebra generated by $X_\nu$ where $\nu\in\Words(A).B$ (point means concatenation here).
Then $\theta_{A,B}|_{\mathfrak h_{\deg_B>0}}$ can inverted back by the prescription
 $X_{\alpha_1\ldots\alpha_n\beta}\mapsto\bo[X_{\alpha_1},\bo\ldots\bo[X_{\alpha_n},X_\beta\bo]\ldots\bo]\bo]$,
 which can be taken as the range is a free Lie algebra.
Ultimately,  this is elimination theory in its simplest form,
 identifying $\mathfrak h_{\deg_B>0}$ to $\mathrm F^{\Lie}_K[X_\nu :\nu\in\Words(A).B]$,
 eliminating the variables $X_\alpha$ $(\alpha\in A)$.
In fact, this above is the Lazard--Shirshov elimination process,
 see \v{S}ir\v{s}ov \cite{S0}, Lazard \cite{L2}, Bourbaki \cite{BX}, Reutenauer \cite{R}.
(It is most popular when $|A|=1$.)
Note that the elimination map $\theta_{A,B}|_{\mathfrak h_{\deg_B>0}}$, in a given multigrade,
 decreases the joint degree in $X$ by $\deg_A$ but it takes the multidegree to a direct sum of full multidegrees
(one has to account the number of letters from $A\cup B$ in various distributions).
This allows to prove the (strong) MW theorems by induction:
\snewpage

\begin{proof}[Proof of the strong MW'  theorem via Lie algebraic elimination.]
We use induction on the joint degrees of the multidegrees.
Let us consider a multidegree.
We can assume that the number of its variables is at least two.
Let us eliminate a variable.
Then the multidegree is isomorphic to the direct sum of some multidegrees of less joint degree.
By induction, this implies freeness in the multidegree immediately.
By induction, in the image monomial bases can be chosen.
Back substitution in the elimination (which inverts the elimination map) takes Lie monomials into Lie monomials.
\end{proof}

\begin{proof}[Proof of the MW  theorem via Lie algebraic elimination.]
It is sufficient to prove that if Lie polynomial $P^{\Lie}(X_1,\ldots,X_n)\neq 0$ is  in a given multidegree,
 then adjoining an external variable $Y$, one has $\bo[P^{\Lie}(X_1,\ldots,X_n),Y\bo]\neq 0$.
(In that case $P^{\ass}(\ad X_1,\ldots,\ad X_n)Y=\bo[P^{\Lie}(X_1,\ldots,X_n),Y\bo]\neq 0$
 provides a good associative representation.)
We use induction on the joint degrees in $X$ of the multidegrees.
We can assume that the variables in $P^{\Lie}$ are of multiplicity at least $1$, and $n\geq2$
 (noncommutative Lie algebras over any $K$ exist after all).
Eliminating any variable $X_i$ does not change the overall shape of the actors but decreases the
 joint degree in $X$, thus, by induction, the statement follows.
\end{proof}

\begin{proof}[Proof of the strong MW''  theorem via Lie algebraic elimination.]
Here we can consider Lie polynomials $Q^{\Lie}(X_1,\ldots,X_n,Y)\neq 0$
 in multidegrees such that every variable has multidegree at least $1$, but in case of $Y$,
 it is exactly $1$. let us call such a Lie polynomial `old' if it has shape
 $\bo[P^{\Lie}(X_1,\ldots,X_n),Y\bo]$.
It is sufficient to prove that, in the given multidegree, the subspace of `old' Lie polynomials
 can be complemented by a subspace of ``new'' Lie polynomials for which there is a basis of shape
 $ \bo[M_1,\ldots,M_k,Y\bo]_{\mathrm L}$ where the $M_i$ are Lie polynomials of the $X_i$.
Again, we prove this by induction on the joint degree in $X$.
We can assume that $n\geq2$, otherwise the statement is simple.
Let us eliminate a variable $X_i$.
In the result the multidegree goes to a direct sum of full multidegree, whose variables
 are some $X_\lambda$ (but at least one) with various multiplicities, and a variable $Y_{i\ldots i}$
 (with possibly empty index) with multiplicity exactly $1$.
One can see, however, that the elimination map (and the inverse back substitution map)
 is set up so that `old' Lie polynomials in the multidegree go exactly to
 `old' Lie polynomials in multidegree involving the variable $Y$ (with empty index).
Thus the complement can be constructed from the complements in the image multidegrees
 involving $Y$ and the full multidegrees involving the  $Y_{i\ldots i}$ with non-empty indices,
 and by back substitution.
In the  multidegrees involving the  $Y_{i\ldots i}$ monomial bases can be chosen, and
 back substitution keeps the bases in the required form.
\end{proof}
 \renewcommand{\qedsymbol}{$\triangle$}
\begin{proof}[Note]
Here in all proofs we proceeded by eliminating only one variable.
\end{proof}
\renewcommand{\qedsymbol}{$\Box$}

Returning to elimination as such, from the Magnus--Witt theorem it also follows that not only
 $\theta_{A,B}|_{\mathfrak h_{\deg_B>0}}$ is an isomorphism but $\theta_{A,B}$ is injective.
This yields the natural isomorphism
 $\mathrm F^{\Lie}_K[X_\xi :\xi\in A\cup B]\simeq \mathrm F^{\cmr}_K[D_\alpha :\alpha\in A]\ltimes
 \mathrm F^{\Lie}_K[X_\nu :\nu\in\Words(A).B]$,
 valid even if $A$ or $B$ is empty.
Several variants are possible.

What made difference here is the grading structure and that $K$-derivations of free Lie $K$-algebras can be
 prescribed arbitrarily in the generating sets.
These are not difficult ideas, yet they require a somewhat proficient abstraction in universal algebraic ideas.
However, if it comes to elimination, then it is also simple to use noncommutative polynomials:

\snewpage
\section{About free Lie algebras. Version 2 (``Classical'' version)}\plabel{sec:free2}
\textbf{Elimination by polynomials.}
Consider the noncommutative polynomial  algebra $\mathrm F_K[E_1,\ldots, E_n, X_1,\ldots,X_m ]$.
Let $\theta$ be the operation  which sends the monomial
\begin{equation}
\underbrace{E_{i_{1,1}}\ldots E_{i_{1,p_1}}X_{j_1}\ldots E_{i_{s,1}}\ldots E_{i_{s,p_s}}X_{j_s}}_{\text{$s$ many occurrences of $X$}} E_{i_{s+1,1}}\ldots E_{i_{s+1,p_{s+1}}}\plabel{eq:word}
\end{equation}
into the polynomial
\[[E_{i_{1,1}},\ldots ,E_{i_{1,p_1}},X_{j_1}]_{\mathrm L}\ldots [E_{i_{s,1}},\ldots ,E_{i_{s,p_s}},X_{j_s}]_{\mathrm L} E_{i_{s+1,1}}\ldots E_{i_{s+1,p_{s+1}}}.\]
\begin{lemma}\plabel{lem:prerec}
The map $\theta$  leaves the multigrading of $\mathrm F_K[E_1\ldots E_n, X_1,\ldots,X_m ]$ invariant.
It acts as an isomorphism in every multigrade.
\begin{proof}
It is obvious that the multigrading is left invariant.
Let us choose an ordering on the letter $E_1,\ldots, E_n, X_1,\ldots,X_m$
 but so that the $E_1,\ldots, E_n$ are all greater than $X_1,\ldots,X_m$.
Then $\theta$ changes all associative monomials by antilexicographically greater associative monomials, thus
 with respect to antilexicographical ordering of the associative monomials
 action of $\theta$ is triangular with $1$'s in the diagonal, thus it is an isomorphism.
\end{proof}
\end{lemma}

A  divided noncommutative polynomial $P(Y_1,\ldots,Y_r| E_1,\ldots, E_n)$ is
 a noncommutative polynomial   where the monomials are of shape $Y_{k_1}\ldots Y_{k_p} E_{i_1}\ldots E_{i_q}$.
\begin{lemma}\plabel{lem:rec}
Suppose that $P(Y_1,\ldots,Y_r\mid E_1,\ldots E_n)$ is a noncommutative polynomial over $K$,
 and assume that the noncommutative polynomial
\[P([E_{i_{1,1}},\ldots ,E_{i_{1,p_1}},X_{j_1}]_{\mathrm L},\ldots,  [E_{i_{r,1}},\ldots ,E_{i_{r,p_r}},X_{j_r}]_\mathrm L
 \mid E_1,\ldots, E_n   )=0,\]
 yet the monomials
\[E_{i_{1,1}}\ldots E_{i_{1,p_1}}X_{j_1}\quad,\quad\ldots\quad,\quad E_{i_{r,1}}\ldots E_{i_{r,p_r}}X_{j_r}\]
 are different from each other.
Then
\[P(Y_1,\ldots,Y_r\mid E_1,\ldots E_n)=0.\]
\begin{proof}
Let us apply the isomorphism $\theta^{-1}$.
This gives
\[P(E_{i_{1,1}}\ldots E_{i_{1,p_1}}X_{j_1},\ldots,  E_{i_{r,1}}\ldots E_{i_{r,p_r}}X_{j_r} \mid E_1,\ldots, E_n   )=0.\]
But then back recognition on the monomials implies $P(Y_1,\ldots,Y_r\mid E_1,\ldots E_n)=0$.
\end{proof}
\end{lemma}
\begin{cor}\plabel{cor:rec}
Suppose that $P(Y_1,\ldots,Y_r)$ is a noncommutative polynomial over $K$,
 and assume that the noncommutative polynomial
\[P([E_{i_{1,1}},\ldots ,E_{i_{1,p_1}},X_{j_1}]_{\mathrm L},\ldots,  [E_{i_{r,1}},\ldots ,E_{i_{r,p_r}},X_{j_r}]_\mathrm L   )=0,\]
 yet the monomials
\[E_{i_{1,1}}\ldots E_{i_{1,p_1}}X_{j_1}\quad,\quad\ldots\quad,\quad E_{i_{r,1}}\ldots E_{i_{r,p_r}}X_{j_r}\]
 are different from each other.
Then
\[P(Y_1,\ldots,Y_r)=0.\eqed\]
\end{cor}

\begin{proof}[Proof of the Magnus--Witt theorem by noncommutative polynomial elimination.]
We have to prove that if $P^{\Lie}\in\mathrm F^{\Lie}_K[X_\lambda:\lambda\in\Lambda]$ evaluates to $0$
 in the commutator expansion in $\mathrm F_K[X_\lambda:\lambda\in\Lambda]$, then
 $P^{\Lie}$ simplifies to $0$ in $\mathrm F^{\Lie}_K[X_\lambda:\lambda\in\Lambda]$.
We can assume that $P^{\Lie}$ is expanded to Lie-monomials, thus it is represented by a non-associative polynomial $P^{\nonass}$.
In that viewpoint, we have to prove that if $P^{\nonass}\in\mathrm F^{\nonass}_K[X_1,\ldots,X_n]$
 evaluates to $0$ in the commutator expansion in $\mathrm F_K[X_1,\ldots,X_n]$, then
 $P^{\nonass}$ can be simplified to $0$ using Lie rules.
We prove the statement by induction on the maximal length $\deg(P^{\nonass})$ of the $\bo[,\bo]$-monomials in  $P^{\nonass}$.
If $\deg P^{\nonass}=0$, then the statement is obvious.
Let us gather the terms of $P^{\nonass}$ into groups $P_\omega^{\nonass}$ corresponding to multigrades.
The various $P_\omega^{\nonass}$ expand to different multigrades in  $\mathrm F_K[X_1,\ldots,X_n]$,
 thus the various $P_\omega^{\nonass}$ must also expand to $0$ in $\mathrm F_K[X_1,\ldots,X_n]$ independently.
Thus it is sufficient to consider the cases $P_\omega^{\nonass}$ separately.
We can assume that $P^{\nonass}_\omega$ has monomials with variables $X_1,\ldots, X_n$
 with multiplicities $i_1,\ldots,i_n\geq1$ respectively.
If $n=1$, then the statement is very easy: in the  $i_1=1$ case the commutator expansion is identical,
 in the $i_1>1$ case $P_\omega^{\nonass}$ obviously reduces to $0$ using Lie rules.
So, assume $n\geq2$.
Then by standard Lie rules we can expand $P_\omega^{\nonass}$ to a Lie-polynomial of
 some $\bo[X_{j_1}, \ldots,X_{j_p} ,X_n\bo]_{\mathrm L}$
 ($j_k<n$) but so that formally the multiplicities of the variables remain.
Thus
\begin{multline}
P^{\nonass}_\omega(X_1,\ldots ,X_n)=\\=Q^{\nonass}_\omega( \bo[X_{j_{1,1}}, \ldots,X_{j_{1,p_1}} ,X_n\bo]_{\mathrm L}, \ldots , \bo[X_{j_{r,1}}, \ldots,X_{j_{r,p_r}} ,X_n\bo]_{\mathrm L} )\mod\Lie,
\plabel{eq:frec}
\end{multline}
 where the sequences $(X_{j_{1,1}}, \ldots,X_{j_{1,p_1}}),\ldots,(X_{j_{s,1}}, \ldots,X_{j_{r,p_r}} ) $ are different from each other,
 while the  multiplicities of  the variables $X_i$ on the two sides are the same.
Nevertheless, the RHS of \eqref{eq:frec} must also expand to $0$ in the commutator evaluation.
But then according to Corollary \ref{cor:rec}, $Q^{\nonass}_\omega(Y_1,\ldots,Y_r)$ also expands to $0$ in the commutator expansion.
Now $\deg Q^{\nonass}_\omega=i_n<\deg P^{\nonass}$ due to the multiplicity structure,
 thus, by induction, we know that  $Q^{\nonass}_\omega(Y_1,\ldots,Y_r)$ expands to $0$ using Lie rules.
But this implies that the RHS of  \eqref{eq:frec} simplifies to $0$ using Lie rules.
So, consequently, also the LHS of  \eqref{eq:frec}.
\end{proof}

Then, by the Magnus--Witt theorem, we see that $\mathrm F^{\Lie}_K[X_1,\ldots,X_n]$
 is multigraded induced from the multigrading in $\mathrm F_K[X_1,\ldots,X_n]$ through commutator evaluation.
(Note that in the previous proof we have not used the gradedness of $\mathrm F^{\Lie}$ but  of $\mathrm F^{\nonass}$.)

\begin{proof}[Proof of the strong Magnus--Witt' theorem  by noncommutative polynomial elimination.]
We can assume that in a multigrade we have variables $X_1,\ldots,X_n$  with multiplicities $i_1,\ldots,i_n\geq1$ respectively.
We proceed by induction on the degree $i_1+\ldots+i_n=i$.
(Due to the previous statement we will use the terms Lie-polynomial and commutator polynomial  synonymously.)
If $i=0$, then the statement is obvious.
Assume that $i\geq1$.
If $n=1$, then the statement is trivial.
So assume $n\geq2$.
Using standard Lie rules, any Lie-polynomial $P$ of multigrade $X_1^{i_1}\ldots X_n^{i_n}$ can be written in form
\[Q( [X_{j_{1,1}}, \ldots,X_{j_{1,p_1}} ,X_n]_{\mathrm L}, \ldots , [X_{j_{r,1}}, \ldots,X_{j_{r,p_r}} ,X_n]_{\mathrm L} )\]
 such that $Q(Y_1,\ldots,Y_r)$ is a Lie-polynomial, and the
 sequences $(X_{j_{k,1}}, \ldots,X_{j_{k,p_k}}) $
 run through every word of length at most $i$ made from $\{X_1,\ldots,X_{n-1}\}$.

Regarding the multigrade structure of $Q$, not every multigrade $Y_1^{j_1}\ldots Y_{r}^{j_r}$ is allowed to appear nontrivially.
But if an multigrade $Y_1^{j_1}\ldots Y_{r}^{j_r}$ is allowed, then every commutator polynomial
 of multigrade $Y_1^{j_1}\ldots Y_{r}^{j_r}$ is allowed to be used in $Q$.
Thus we obtain that $P$ is of shape
\[\sum_{Y_1^{j_1}\ldots Y_{r}^{j_r} \text{is an allowed multigrade}}
Q_{Y_1^{j_1}\ldots Y_{r}^{j_r}}( [X_{j_{1,1}}, \ldots,X_{j_{1,p_1}} ,X_n]_{\mathrm L}, \ldots , [X_{j_{r,1}}, \ldots,X_{j_{r,p_r}} ,X_n]_{\mathrm L} )\]
 such that $Q_{Y_1^{j_1}\ldots Y_{r}^{j_r}}$ is of multigrade $Y_1^{j_1}\ldots Y_{r}^{j_r}$.
On the other hand, this description is unique in terms of the commutator polynomials
 $Q_{Y_1^{j_1}\ldots Y_{r}^{j_r}}$ due to Corollary \ref{cor:rec}.
Hence, the situation decomposes in allowed multigrades.
Thus, in particular, if
\begin{equation}
Q_\lambda(Y_1,\ldots, Y_r)\quad:\quad\lambda\in\Lambda_{Y_1^{j_1}\ldots Y_{r}^{j_r}}
\notag
\end{equation}
form systems of base monomials in the allowed multigrades $Y_1^{j_1}\ldots Y_{r}^{j_r}$, then the elements
\begin{equation}
Q_\lambda( [X_{j_{1,1}}, \ldots,X_{j_{1,p_1}} ,X_n]_{\mathrm L}, \ldots , [X_{j_{r,1}}, \ldots,X_{j_{r,p_r}} ,X_n]_{\mathrm L} )\,:
\,\lambda\in\bigcup_{\substack{Y_1^{j_1}\ldots Y_{r}^{j_r} \\\text{is an allowed}\\\text{ multigrade}}}^\circ \Lambda_{Y_1^{j_1}\ldots Y_{r}^{j_r}}
\notag
\end{equation}
 form a system of base monomials for multigrade  $X_1^{i_1}\ldots X_n^{i_n}$.
However, for allowed multigrades, $\deg Y_1^{j_1}\ldots Y_{r}^{j_r}=i_n<\deg X_1^{i_1}\ldots X_n^{i_n}$,
 so by induction we have monomial bases in allowed multigrades.
It is also clear that this process can be made independent from the actual coefficients $K$.
\end{proof}

\begin{proof}[Proof of the strong Magnus--Witt'' theorem  by noncommutative polynomial elimination.]
By the properties of $\theta$ (Lemma \ref{lem:prerec} and Lemma \ref{lem:rec}),
 any  polynomial $P$ of multigrade $X_1^{i_1}\ldots X_n^{i_n}$ $(n\geq2)$ can be written in form
\[Q( [X_{j_{1,1}}, \ldots,X_{j_{1,p_1}} ,X_n]_{\mathrm L}, \ldots , [X_{j_{r,1}}, \ldots,X_{j_{r,p_r}} ,X_n]_{\mathrm L}
\mid X_1,\ldots,X_{n-1} )\]
 such that $Q(Y_1,\ldots,Y_r\mid X_1,\ldots,X_{n-1} )$ is a divided polynomial, and such that the
 sequences $(X_{j_{k,1}}, \ldots,X_{j_{k,p_k}}) $
 run through every word of length at most $i$ made from $\{X_1,\ldots,X_{n-1}\}$.
The allowed polynomials $Q$ are direct  sums of full admissible multidegrees of divided polynomials.
From the uniqueness of $Q$ we know that if $P$ is a commutator polynomial,
 then in $Q(Y_1,\ldots,Y_r\mid X_1,\ldots,X_{n-1} )$ the $X$-degree is $0$.
Then the complement can be constructed from the admissible multidegrees with $\deg_X>0$, and
 from the complements in admissible multidegrees with $\deg_X=0$, by induction, by back substituting $Y_k$.
In the  admissible multidegrees with $\deg_X>0$ we can choose an associative monomial basis,
 and the back substitution process is compatible to shape of bases indicated.
\end{proof}
 \renewcommand{\qedsymbol}{$\triangle$}
\begin{proof}[Note]
Here in all proofs we proceeded by eliminating all but one variable.
But we could have done otherwise easily.
Indeed, elimination by noncommutative polynomials is due to Magnus \cite{MM}, and originally, he eliminates only one variable.
(Witt \cite{W} uses the PBW theorem for free Lie algebras over fields first; and then,
the treatment of the integral case is more in group theoretic terms.)
It must be clear that elimination by derivations and elimination by noncommutative polynomials
 are two sides of the same coin, but, written down, they are somewhat different.
For our purposes, elimination by noncommutative polynomials is very natural.
In terms of generalizations, however, derivations provide a more extendible approach,
cf. Bremner, Dotsenko \cite{BD}.
\end{proof}
\renewcommand{\qedsymbol}{$\Box$}
\snewpage

\section{$\mathrm F_{K}^{\Lie}$ via $\mathrm F_{\mathbb Q}^{\Lie}$}\plabel{sec:via}
Firstly, we give a direct proof of the PBW theorem for  free Lie algebras over $\mathbb Q$.
(The argument works for any field of characteristic $0$.)
We will use the fact that free Lie algebras are multigraded. The proof will be sketchy as we rely on familiar arguments.

We define a  Lie-permutation $Ii$ of $\{1,\ldots,n\}$ as the following data.
It is a partition   $I_1\dot\cup\ldots \dot\cup I_s=\{1,\ldots,n\}$ such that $\max I_1<\ldots<\max I_s$,
 and finite sequences $i_{k,1},\ldots,i_{k,p_k}$ such that $\{i_{k,1},\ldots,i_{k,p_k}\}=I_k$, $p_k=|I_k|$ and $i_{k,p_k}=\max I_k$.
Let $\mathbb QX_{\Sigma_n}$ be the vector space spanned by the noncommutative monomials
 $X_{\sigma(1)}\ldots X_{\sigma(n)}$ in the corresponding noncommutative polynomial ring over $\mathbb Q$.
From \cite{LLX}, we quote the easy
\begin{prop}\plabel{prop:deco}
Any element of $\mathbb QX_{\Sigma_n}$ can uniquely be written in the form
\begin{equation}
\sum_{Ii\text{ is a Lie-permutation of }\{1,\ldots,n\}} a_{Ii}
[X_{i_{1,1}},\ldots,X_{i_{1,p_1}}]_{\mathrm L}\cdot_\Sigma\ldots \cdot_\Sigma [X_{i_{s,1}},\ldots,X_{i_{s,p_s}}]_{\mathrm L}
\plabel{eq:decoe}
\end{equation}
 where $a_{Ii}\in\mathbb Q$.
(Here we used ordinary commutators and symmetrized products.)
\qed
\end{prop}
\begin{prop} The PBW theorem holds for  $\mathfrak g=\mathrm F_{\mathbb Q}^{\Lie}[X_\lambda:\lambda\in\Lambda]$.
\begin{proof}[Sketch of proof]
We will prove the symmetric global formulation.
Consider
\[\boldsymbol m_\Sigma:\textstyle{\bigotimes_\Sigma\mathfrak g}\rightarrow\mathcal U\mathfrak g\simeq
\mathrm F_{\mathbb Q}[X_\lambda:\lambda\in\Lambda]. \]
Both sides are naturally multigraded, and the map is compatible with them.
Thus, it is sufficient to prove isomorphism (that is injectivity) between them in every multigrade separately.
If in the multigrade every variable has multiplicity at most one, then the injectivity holds due to Proposition \ref{prop:deco}.
Regarding higher multigrades, assume that in  multigrade $X_1^{i_1}\ldots X_s^{i_s}$,
\[P(\underbrace{X_1,\ldots X_1}_{i_1\text{ many terms}},\ldots,\underbrace{X_s,\ldots X_s}_{i_s\text{ many terms}})\in \textstyle{\bigotimes_\Sigma\mathfrak g}\]
 is such that it is not zero but evaluates to zero in $ \mathrm F_{\mathbb Q}[X_\lambda:\lambda\in\Lambda]$.
Here $P$ was written such that every variable appears according to the multiplicity (for a  monomial decomposition).
Then the polarization
\begin{equation}
\frac1{n!}\sum_{\sigma\in \Sigma_n}P(\underbrace{X_{1,\sigma(1)},\ldots X_{1,\sigma(i_1)}}_{i_1\text{ many terms}},\ldots,\underbrace{X_{s,\sigma(n-i_s+1)},\ldots X_{s,\sigma(n)}}_{i_s\text{ many terms}})
\plabel{eq:polar}
\end{equation}
($n=i_1+\ldots+i_s$) is also not zero, because its depolarization is nonzero.
On the other hand, it evaluates to
 the polarization of $0$, i.~e.~to $0$ as a noncommutative polynomial.
This contradicts to the injectivity of the multigrades without multiplicity $\geq2$.
\end{proof}
\end{prop}

The PBW theorem for $\mathfrak g=\mathrm F^{\Lie}_{\mathbb Q}[X_1,\ldots, X_n]$ yields $\mu^{\Lie}_n(X_1,\ldots,X_n)$
 as the component of degree 1 of $(\bo m_\Sigma )^{-1}(X_1\ldots X_n)$, that is, as a canonical projection, cf. \cite{LLX}.
Using these a formula can be given for $(\bo m_\Sigma )^{-1}$ in the general case $\mathbb Q\subset K$.
Altogether, this is not so different from or any simpler than the universal algebraic (second) proof in \cite{LLX};
 but both observations that PBW theorem for $\mathrm F^{\Lie}_{\mathbb Q}$ is somewhat easier to prove
 and that the PBW theorem for $\mathrm F^{\Lie}_{\mathbb Q}$ implies the symmetric PBW theorem
 are sort of notable.
\snewpage

We will use Proposition \ref{prop:CMag}, the strong Magnus--Witt' theorem, in order to prove
\begin{prop} The PBW theorem holds for  $\mathfrak g=\mathrm F_{K}^{\Lie}[X_\lambda:\lambda\in\Lambda]$.
\begin{proof}First, let us consider the case $K=\mathbb Z$.
We know that
$\mathrm F_{\mathbb Z}^{\Lie}[X_\lambda:\lambda\in\Lambda]$ is a free $\mathbb Z$-module,
 and  $ \mathrm F_{\mathbb Z}^{\Lie}[X_\lambda:\lambda\in\Lambda]\otimes\mathbb Q
 \simeq\mathrm F_{\mathbb Q}^{\Lie}[X_\lambda:\lambda\in\Lambda]$
 naturally.
Thus, starting from a basis of $\mathrm F_{\mathbb Z}^{\Lie}[X_\lambda:\lambda\in\Lambda]$, we see that
 $\bigotimes_\leq \mathrm F_{\mathbb Z}^{\Lie}[X_\lambda:\lambda\in\Lambda]$ embeds to
 $\bigotimes_\leq \mathrm F_{\mathbb Q}^{\Lie}[X_\lambda:\lambda\in\Lambda]$.
The latter one evaluates in $\mathcal U \mathrm F_{\mathbb Q}^{\Lie}[X_\lambda:\lambda\in\Lambda]$ injectively, thus
 $\bigotimes_\leq \mathrm F_{\mathbb Z}^{\Lie}[X_\lambda:\lambda\in\Lambda]$ evaluates
 (taking Lie-commutator into commutator, tensor product into ordinary product) in some ring injectively.
This implies that $\bigotimes_\leq \mathrm F_{\mathbb Z}^{\Lie}[X_\lambda:\lambda\in\Lambda]$
 must evaluate in the universal enveloping  $\mathcal U \mathrm F_{\mathbb Z}^{\Lie}[X_\lambda:\lambda\in\Lambda]$ algebra injectively.

Let us now consider the general case.
We know that the evaluation yields an isomorphism
 $\bigotimes_\leq \mathrm F_{\mathbb Z}^{\Lie}[X_\lambda:\lambda\in\Lambda]\simeq \mathrm F_{\mathbb Z}[X_\lambda:\lambda\in\Lambda]$
 of free $\mathbb Z$-modules.
But then evaluation (i.~e.~the process which sends Lie-commutators to commutators and tensor product to products)
 gives an isomorphism $K\otimes\bigotimes_\leq
 \mathrm F_{\mathbb Z}^{\Lie}[X_\lambda:\lambda\in\Lambda]\simeq \mathrm F_{K}[X_\lambda:\lambda\in\Lambda]$.
On the other hand, $K\otimes\bigotimes_\leq \mathrm F_{\mathbb Z}^{\Lie}[X_\lambda:\lambda\in\Lambda]
 \simeq \bigotimes_\leq \mathrm F_{K}^{\Lie}[X_\lambda:\lambda\in\Lambda]$
 naturally, and compatibly with evaluation. That proves that $\bigotimes_\leq
 \mathrm F_{K}^{\Lie}[X_\lambda:\lambda\in\Lambda]$ evaluates to $\mathrm F_{K}[X_\lambda:\lambda\in\Lambda]$
 isomorphically.
\end{proof}
\end{prop}
The bottleneck  in the previous line of arguments was that we could establish  only the case of $\mathbb Q$ first.
This is because we used symmetrization.
In fact, symmetrization was used in two instances.
The first one was the use of Proposition \ref{prop:deco}.
However, this proposition can be modified easily:
Let $KX_{\Sigma_n}$ be the $K$-module spanned by the noncommutative monomials
$X_{\sigma(1)}\ldots X_{\sigma(n)}$ in the corresponding noncommutative polynomial ring over $  K$.
\begin{prop}\plabel{prop:decovar}
Any element of $KX_{\Sigma_n}$ can uniquely be written in the form
\[
\sum_{Ii\text{ is a Lie-permutation of }\{1,\ldots,n\}} a_{Ii}
[X_{i_{s,1}},\ldots,X_{i_{s,p_s}}]_{\mathrm L}\cdot \ldots \cdot[X_{i_{1,1}},\ldots,X_{i_{1,p_1}}]_{\mathrm L}
\]
 where $a_{Ii}\in K$.
(Here we used ordinary commutators.)
\begin{proof}[Sketch of proof]
It can be seen that the antilexicographically highest monomial term in
$[X_{i_{s,1}},\ldots,X_{i_{s,p_s}}]_{\mathrm L}\cdot \ldots \cdot[X_{i_{1,1}},\ldots,X_{i_{1,p_1}}]_{\mathrm L}$
 is just $(X_{i_{s,1}}\cdot\ldots\cdot X_{i_{s,p_s}})\cdot \ldots \cdot(X_{i_{1,1}}\cdot\ldots\cdot X_{i_{1,p_1}})$,
 and it is with coefficient $1$.
This implies that we have a triangular change a basis relative to the ordinary monomial basis.
\end{proof}
\end{prop}
The second instance to use symmetrization was to deal with the case of multiplicities in the variables.
That cannot be done particularly easily (although the strong Magnus--Witt' theorem is not that difficult).
Two possible strategies are as follows:

One strategy is eliminating the rational case by using the Magnus--Witt theorems more effectively.
This will be done in Section \ref{sec:freeMW}.

Another strategy is to make Proposition \ref{prop:decovar} work for multiplicities by better accounting
 (thus recreating the Magnus--Witt theorems in more precise form).
This will be realized in Section \ref{sec:freedirect}. (Actually, we change the ordering of variables,
 and instead of `antilexicographically highest' we will have `antilexicographically smallest',
 but such changes [like $X_{n}\leftrightarrow X_{-n}$] can be considered as [inconveniently] usual.)
\snewpage
\section{The $\mathrm F_{K}^{\Lie}$ case from the MW theorems and the coproduct structure}\plabel{sec:freeMW}
Here the idea is to use the coproduct structure on the $\mathcal U\mathfrak g$.
In order to deal with case $\mathbb Q \not\subset K$, we use the collective symmetric coproduct.
For pedagogical reason, first we define it for $\mathrm F_{K}[X_\lambda:\lambda\in\Lambda]$ such that
\[\bo\Delta:\mathrm F_{K}[X_\lambda:\lambda\in\Lambda]\rightarrow\bigoplus_{n=0}^\infty\textstyle{\bigodot^n}
\underbrace{\mathrm F_{K}^{\geq 1}[X_\lambda:\lambda\in\Lambda]}_{\text{no constant terms}}.\]
 by the recipe
\[\bo\Delta
 (X_{\lambda_1} \ldots X_{\lambda_n})=\sum_{\substack{I_1\dot\cup\ldots \dot\cup I_s=\{1,\ldots,n\}\\I_k=\{i_{k,1},\ldots,i_{k,p_k}\}\neq
\emptyset \\i_{k,1}<\ldots<i_{k,p_k}}} \,
  (X_{\lambda_{i_{1,1}}}\ldots X_{\lambda_{i_{1,p_1}}})\odot\ldots\odot
  (X_{\lambda_{i_{s,1}}},\ldots ,X_{\lambda_{i_{s,p_s}}}).
\]
This behaves well with respect relabeling and substitution of variables.
It also behaves well  with respect substitution of commutators:
\[\bo\Delta(X_{\lambda_1} \ldots X_{\lambda_{k-1}}[X_{\lambda_{k,1}},X_{\lambda_{k,2}} ] X_{\lambda_{k+1}} \ldots X_{\lambda_n})
=\subs_{Y=[X_{\lambda_{k,1}},X_{\lambda_{k,2}} ]}
\bo\Delta(X_{\lambda_1} \ldots X_{\lambda_{k-1}}Y X_{\lambda_{k+1}} \ldots X_{\lambda_n}) \]
 (understood appropriately).
Indeed, one has to consider the distribution of $X_{\lambda_{k,1}},X_{\lambda_{k,2}}$,
 and one can deduce that the have to fall into the same ``slots'' in order to avoid canceling.

The construction also works for $\mathrm F_K[\mathfrak g]=\textstyle{\bigotimes^n\mathfrak g}$,
 and we find (with some abuse of notation)
\begin{multline*}\bo\Delta_\otimes(Z_{ 1}\otimes \ldots\otimes Z_{ {k-1}}\otimes[Z_{ {k,1}},Z_{ {k,2}} ]_{\otimes}\otimes Z_{ {k+1}}\otimes \ldots \otimes Z_{ n})
=\\=\subs_{Y=[Z_{ {k,1}},Z_{ {k,2}} ]_\otimes}
\bo\Delta_\otimes(Z_{ 1}\otimes \ldots\otimes Z_{ {k-1}}\otimes Y \otimes Z_{ {k+1}}\otimes \ldots  \otimes Z_{ n}).
\end{multline*}
Therefore,
\begin{multline*}\bo\Delta_\otimes(Z_{ 1}\otimes \ldots\otimes Z_{ {k-1}}\otimes
\left([Z_{ {k,1}},Z_{ {k,2}} ]_{\otimes}-\bo[Z_{k,1},Z_{k,2}\bo]\right)\otimes Z_{ {k+1}}\otimes \ldots \otimes Z_{ n})
=\\=\subs_{Y=[Z_{ {k,1}},Z_{ {k,2}} ]_\otimes-\bo[X_{k,1},X_{k,2}\bo]}
\bo\Delta_\otimes(Z_{ 1}\otimes \ldots\otimes Z_{ {k-1}}\otimes Y \otimes Z_{ {k+1}}\otimes \ldots  \otimes Z_{ n}).
\end{multline*}
This means that the symmetric coproduct is compatible to the factorization structure
 of the universal enveloping algebra, yielding a map
\[\bo\Delta:\mathcal U\mathfrak g\rightarrow\bigoplus_{n=0}^\infty\textstyle{\bigodot^n}
\mathcal U^{\geq1}\mathfrak g.\]
Note compatibility with the filtration structure: $\bo\Delta$ takes
$\mathcal U^{\leq s}\mathfrak g$ to $\bigoplus_{n=0}^s\textstyle{\bigodot^n} \mathcal U^{\geq1}\mathfrak g$.
(This second, more general definition of $\bo\Delta$ will dominate our thing thinking in the following.)

\begin{prop} The PBW theorem holds for  $\mathfrak g=\mathrm F_{K}^{\Lie}[X_\lambda:\lambda\in\Lambda]$.
\begin{proof}
In the case $\mathfrak g=\mathrm F_{K}^{\Lie}[X_\lambda:\lambda\in\Lambda]$, by the Magnus--Witt' theorem, we know that
 it is a free Lie $K$-module, thus the global setup of the basic PBW theorem applies.
We can consider the composition of maps
\[\mathrm F_{K}^{\geq 1}[X_\lambda:\lambda\in\Lambda]\xrightarrow{\pi}
\mathrm F_{K}^{\cmr}[X_\lambda:\lambda\in\Lambda]\xrightarrow{\varkappa}
\mathrm F_{K}^{\Lie}[X_\lambda:\lambda\in\Lambda],\]
 where $\pi$ is (the restriction of) a projection map from the Magnus--Witt'' theorem (multi\-gra\-de-wise, in particular, killing $1$),
 and $\varkappa$ is the inverse isomorphism from the Magnus--Witt theorem.
Thus $\varkappa\circ\pi$ is a splitting map for the commutator evaluation
 $\iota:\mathrm F_{K}^{\Lie}[X_\lambda:\lambda\in\Lambda]\rightarrow \mathrm F_{K} [X_\lambda:\lambda\in\Lambda]$.
This yields a map
\[\bo\Delta_\pi=\varkappa_*\circ \pi_*\circ\bo\Delta:\mathcal U\mathfrak g\rightarrow\bigoplus_{n=0}^\infty
\textstyle{\bigodot^n}\mathfrak g\equiv \textstyle{\bigodot}\mathfrak g.\]
 (both projection to a direct component and isomorphisms are compatible to $\odot$).
In order to show that $\bo m_\leq$ is injective, it is sufficient to show that $\bo\Delta_\pi\circ\bo m_\leq$ is injective.
Assume that $S\in \textstyle{\bigotimes_\leq\mathfrak g}$ such that its highest nontrivial degree is $n$.
But then, by the filtration property,
 the projection of $\bo\Delta_\pi\circ\bo m_\leq(S)$ to component $\textstyle{\bigodot^n}\mathfrak g  $
 is the same as the projection of $S$ to the component $\textstyle{\bigotimes^n_\leq}\mathfrak g$
 and the natural isomorphism $\textstyle{\bigotimes^n_\leq}\mathfrak g\rightarrow \textstyle{\bigodot^n}\mathfrak g $ taken.
This proves the injectivity of $\bo\Delta_\pi\circ\bo m_\leq$.
\end{proof}
\end{prop}
The proof simplifies if $K$ is a field but, even so, it requires the MW theorem
(which is relatively effortless for $\mathbb Q\subset K$).
In the view of the PBW theorem, it is easy to see
\begin{cor}
If $S\in \mathrm F_{K}[X_\lambda:\lambda\in\Lambda]$ is without a constant term, then
 $S$ is a commutator polynomial if an only if $\bo\Delta(S)=S$.
\qed
\end{cor}
The latter condition is checkable but computationally demanding.  (Cf. Cohn \cite{CC2}.)
\snewpage
\section{The $\mathrm F_{K}^{\Lie}$ case directly}\plabel{sec:freedirect}
\textbf{Free PBW word bases.}
We say that a free PBW word basis is the following data.
We will consider  words formed from an alphabet $\Lambda$.
\begin{enumerate}
\item[(A1)] Some words  should be called primitive.
\item[(A2)] To any primitive word $w$ a $\bo[,\bo]$-monomial $P^{\Lie}_w$ should be associated such that the variables
 $X_\lambda$ of $P^{\Lie}_w$ correspond to the $\lambda$ in $w$ with the same multiplicity.
\item[(A3)] The Lie-polynomials $P^{\Lie}_w$ should generate $\mathrm F_{K}^{\Lie}[X_\lambda:\lambda\in\Lambda]$ as a $K$-module.
\item[(A4)] The primitive words should be endowed by an ordering $\preccurlyeq$ such that
 every word $w$ uniquely decomposes to a concatenation of primitive words $w=w_1\ldots w_s$ such that
 $w_1  \preccurlyeq w_2 \preccurlyeq\ldots\preccurlyeq w_s$.
\item[(A5)] To any word decomposed as above we associate the noncommutative polynomial
\[P^\ass_w=P^\ass_{w_1}\ldots P^\ass_{w_s},\]
 where $P^\ass_{w_i}$ is the commutator evaluation $P^{\Lie}_{w_i}$.
\item[(A6)] The noncommutative polynomials $P^\ass_w$ should be independent in $\mathrm F_{K}[X_\lambda:\lambda\in\Lambda]$.
\end{enumerate}

Our first observation is that the existence of a free PBW word basis implies the basic PBW theorem for
 $\mathrm F_{K}^{\Lie}[X_\lambda:\lambda\in\Lambda]$.
Indeed, considering (A3) and (A6), we see that that  $P_w^{\Lie}$ should be a basis of $\mathrm F_{K}^{\Lie}[X_\lambda:\lambda\in\Lambda]$.
Then, due to (A3), every element in $\mathcal U \mathrm F_{K}^{\Lie}[X_\lambda:\lambda\in\Lambda]$
 can be brought into a combination of products
 $P_{w_1}^{\Lie},\ldots,P_{w_s}^{\Lie}$ such that $w_1\preccurlyeq w_2\preccurlyeq\ldots\preccurlyeq w_s$,
 by the usual basic rearrangement process.
Such products are then independent in the universal enveloping algebra,
 as they are independent in the noncommutative polynomial evaluation, due to (A6).
This establishes the global form of the  basic PBW theorem with respect to a specific ordering.
But then the local form holds, which implies the global form in general.

\textbf{A variant of Shirsov's construction.}
Next, we will find a free PBW word basis.
There are plenty of such bases but their combinatorics might be confusing for the first sight.
Thus we restrict to a relatively simple construction.
Let $\leq_\Lambda$ be any ordering on $\Lambda$.

Firstly, we need a \textit{breaking pattern}.
We break a finite $\Lambda$-word $w$ with respect to $\leq_\Lambda$ as follows.
The definition is recursive with respect to the length of the words.
We identify the smallest letter $\alpha$ in $w$.
Then $w$ reads as
\begin{equation}
w=i_{1,1},\ldots, i_{1,p_1},\alpha,\ldots, i_{s,1},\ldots ,i_{s,p_s}, \alpha, i_{s+1,1},\ldots, i_{s+1,p_{s+1}}
\end{equation}
 ($s$ many occurrences of $\alpha$).
If the word contains only $\alpha$, then we break the word to letters completely.
If not, then we surely break after the last occurrence after of $\alpha$, and we might break ($\mid_{?}$)
 after other occurrences of $\alpha$, and we might break after the last occurrence of $\alpha$ ($\mid_{??}$):
\begin{equation}
 i_{1,1},\ldots, i_{1,p_1},\alpha\biggr|_?\ldots\biggr|_? i_{s-1,1},\ldots ,i_{s-1,p_{s-1}},
 \alpha\biggr|_?i_{s,1},\ldots ,i_{s,p_s}, \alpha\Biggr| i_{s+1,1}\biggr|_{??}i_{s+1,2}\biggr|_{??}\ldots \biggr|_{??}i_{s+1,p_{s+1}}.
 \notag
\end{equation}
Regarding  $\mid_{??}$ the rule is simple: we break as we would break $i_{s+1,1},\ldots, i_{s+1,p_{s+1}}$.
Regarding  $\mid_?$ the rule is more difficult: Consider the word of tuples
\begin{equation}
w/\alpha=(i_{1,1},\ldots, i_{1,p_1},\alpha) \ldots  (i_{s-1,1},\ldots ,i_{s-1,p_{s-1}},\alpha) (i_{s,1},\ldots ,i_{s,p_s},\alpha).
\plabel{s:elso}\end{equation}
The letters of this word are all elements of the $\Lambda/\alpha$,
 the set of finite tuples from $\Lambda$, whose single minimal element is $\alpha$, and it is the last one.
We say that \eqref{s:elso} is the condensation of the first part of $w$.
We  choose the ordering that $\leq_{\Lambda/\alpha}$
 such that it is the restriction of the $\leq_\Lambda$-antilexicographic ordering.
Then the breaking places of $\mid_?$ are defined to correspond to the breaking places of
 \eqref{s:elso} with respect to $\leq_{\Lambda/\alpha}$.
One can prove by induction on word length that the breaking mechanism is well-defined, as the condensed words always get shorter.
A word is primitive if it does not break (so it gets condensed to a single letter without ever breaking,
 in particular, the latest letter is the maximal).
By that we have (A1).
We call this the canonical decomposition.
Then there is a natural ordering $\preccurlyeq_\Lambda$ defined between primitive $\Lambda$-words as follows:
 $w_1\preccurlyeq_\Lambda w_2$ is the smallest/last letter of $w_1$ is smaller than  the smallest/last letter of $w_2$;
 or if the last letters $\alpha$ are equal and $w_1/\alpha\preccurlyeq_{\Lambda/\alpha} w_2/\alpha$.
By induction on word length, one can prove that the breaking decomposition is non-strictly
 $\preccurlyeq_\Lambda$-increasing, so the  existence part of (A4) are established
 with respect to $ \preccurlyeq_\Lambda$as $\preccurlyeq $ by the canonical decomposition.
The unicity part of (A4) also follows from induction on word length:
If $w=\tilde w_1\ldots\tilde  w_s$ such that
 $\tilde w_1  \preccurlyeq_\Lambda \tilde w_2 \preccurlyeq_\Lambda\ldots\preccurlyeq_\Lambda \tilde w_s$
 is another decomposition the smallest/last letters of the words are in relation
 $\tilde\alpha_1\leq_\Lambda\tilde\alpha_2\leq_\Lambda \ldots\leq_\Lambda\tilde\alpha_s$.
Thus, there is a break after $r$ where $r$ is the biggest index such that $\tilde\alpha_1=\tilde \alpha_r$.
If $r<s$, then the induction hypothesis applies.
If $\alpha\equiv \alpha_1=\ldots=\alpha_s$
 then $w$ can break only after the $\alpha$'s, and using the fact, that $\tilde w_i \preccurlyeq_\Lambda \tilde w_j
 \Leftrightarrow \tilde w_i/\alpha \preccurlyeq_{\Lambda/\alpha} \tilde w_j/\alpha$, we see that the decomposition problem of $w$
 is equivalent to the decomposition problem of $w/\alpha$.
If $w$ is   of only $\alpha$, then the decomposition must be complete as the word $\alpha$ is regular.
If $w$ is not of only $\alpha$, then $w/\alpha$ is shorter than $w$, thus the induction hypothesis applies.
This establishes unicity part of (A4).

Now, $\preccurlyeq_\Lambda$ defined for primitive words can be extend to  $\preccurlyeq_\Lambda^\ext$ for all words:
Let $w \preccurlyeq_\Lambda^\ext \tilde w$ if the decomposition of $w$ to primitives is antilexicographically
 smaller or equal than decomposition of $\tilde w$ to primitives.
(I.~e.~$\preccurlyeq_\Lambda$-big primitives decide.)
One can also view $\preccurlyeq_\Lambda^\ext$ in another way:
Let $(\Lambda|\alpha)$ be set of words whose last element is $\alpha$ and minimal (but possibly empty).
Note that for $w,\tilde w\in (\Lambda|\alpha)$ it holds that
 $w \preccurlyeq_\Lambda^\ext\tilde w \Leftrightarrow w/\alpha \preccurlyeq_{\Lambda/\alpha}^\ext\tilde w/\alpha$.
One can immediately see that $(\Lambda|\alpha) $ is in bijection to $\mathrm{Words( \Lambda/\alpha)}$ by concatenation.
Then there is a unique decomposition
\begin{equation*} w= \ldots (w|\beta)\ldots (w|\gamma)\ldots \end{equation*}
 where in the product $\ldots,\beta,\ldots,\gamma,\ldots$ ranges over all $\Lambda$, such that $\beta<\gamma$
 and $(w|\beta)\in (\Lambda|\beta)$, etc.
Of course, all but finitely many components are empty.
Then $(w|\beta)$ is just the product of the primitive components of $w$, whose last letter is $\beta$.
Let $w/\beta$ be the condensation $(w|\beta)$ with respect to $\beta$.
(This is consistent of the notation before.)
Then one can see that $w\prec_\Lambda^\ext \tilde w$ if and only if there is $\beta$ such that
 $w/\beta \prec^\ext_{\Lambda/\beta}\tilde w/\beta$ but $w/\gamma =\tilde w/\gamma$ for all $\gamma>\beta$.
Together with the observation that for pure $\beta$-words length decides, this characterizes
 $\preccurlyeq_\Lambda^\ext$ without referring to primitives.
At this point, using the previous observations, it is easy to see that
 $\preccurlyeq_\Lambda^\ext$ is   the antilexicographic ordering, and therefore $\preccurlyeq_\Lambda$
 is the restriction of the antilexicographic ordering with respect to the set of primitive words.
\snewpage

Secondly, we need an \textit{evaluation pattern}.
Let us consider a primitive word $w$.
Let us consider how it condenses to a single letter.
What is yields is a single element of $\Lambda$, or a tuple of elements of $\Lambda$ or tuples  of  tuples of $\Lambda$ or\ldots.
For example,
\[  w= \lambda_4\lambda_1 \lambda_1 \lambda_2,\lambda_3    \lambda_1 \lambda_3 \lambda_2    \lambda_1 \lambda_1 \rightsquigarrow\]
\[   (\lambda_4,\lambda_1)(\lambda_1) (\lambda_2,\lambda_3,   \lambda_1)(\lambda_3,\lambda_2   , \lambda_1)(\lambda_1) \rightsquigarrow\]
\[  \bigl( (\lambda_4,\lambda_1),(\lambda_1)\bigr) \bigl((\lambda_2,\lambda_3,
 \lambda_1),(\lambda_3,\lambda_2   , \lambda_1),(\lambda_1)\bigr) \rightsquigarrow\]
\[\Bigl(  \bigl( (\lambda_4,\lambda_1),(\lambda_1)\bigr),\bigl((\lambda_2,\lambda_3,
 \lambda_1),(\lambda_3,\lambda_2   , \lambda_1),(\lambda_1)\bigr)\Bigr) .\]
Then evaluation is set as follows: We replace $\lambda\in\Lambda$ by $X_\lambda$,
 and $(*_1,\ldots,*_n)$ by left-iterated commutators $[ *_1,\ldots,*_n]_{\mathrm L}$.
Regarding the previous example, it yields
\[P^{\Lie}_w=\Bigl[  \bigl[ [X_{\lambda_4},X_{\lambda_1}]_{\mathrm L},
X_{\lambda_1} \bigr]_{\mathrm L},\bigl[[X_{\lambda_2},X_{\lambda_3},   X_{\lambda_1}]_{\mathrm L},
[X_{\lambda_3},X_{\lambda_2}   , X_{\lambda_1}]_{\mathrm L}, X_{\lambda_1} \bigr]_{\mathrm L}\Bigr]_{\mathrm L} .\]
This establishes (A2), and extension with (A5) is immediate.

Then the generating statement (A3) follows by induction (on formal multigrade in $\mathrm F_K^{\nonass}$) using the standard fact that
 Lie-monomials containing $X_n$ are Lie-polynomials of $\bo[X_{j_1}, \ldots,X_{j_k} ,X_\alpha\bo]_{\mathrm L}$.
Indeed, replacing $\bo[X_{j_1}, \ldots,X_{j_k} ,X_\alpha\bo]_{\mathrm L}$ by $X_{(j_1,\ldots,j_k,\alpha )}$,
 and applying induction on degree, the statement follows.

$\Lambda$-words correspond to associative monomials by the recipe
 $w=\lambda_1\ldots\lambda_k\mapsto X_w=X_{\lambda_1}\ldots X_{\lambda_k}$.
That induces the corresponding antilexographic ordering $\preccurlyeq^{\mathrm{mon}}$ on the monomials.
We claim that the $\preccurlyeq^{\mathrm{mon}}$-smallest
 monomial in $P^\ass_w$ with nonzero coefficient is $X_w$, and the corresponding coefficient is $1$.
(This is the ``triangular property''.)
Again, this is easy to prove on word length.
If $w$ decomposes to more than one primitive words, then this follows from the properties of the
 antilexicographical ordering (composed on words with fixed lengths).
If $w$ is primitive and not length $1$, then it is easy to see that in
 expanding   $\bo[X_{j_1}, \ldots,X_{j_k} ,X_\alpha\bo]_{\mathrm L}$ with ${j_1}, \ldots, {j_k}>_\Lambda\alpha$
 the possibly optimal minimal choice is $X_{j_1} \ldots,X_{j_k}X_\alpha$ with respect to
 antilexicographical ordering, and it realizes the indicated value
 when inferring from index case $\Lambda/\alpha$ (with smaller word length).
The independence statement (A6) follows immediately by this triangular property.

\textbf{Note.}
Unbreaking under condensation can reformulated so that
 the primitive words are those which are strictly minimal among their cyclic permutations.
Thus they are exactly the antilexicographical Lyndon--Shirshov words.
See Lyndon \cite{Lyn}, \v{S}ir\v{s}ov \cite{S2}, Melan\c{c}on,  Reutenauer  \cite{MR}.
(Reversing orders and interchanging left and right are usual ambiguities in these matters.)
Free PBW bases were first constructed in the context of Hall--Shirshov bases (cf.~the Introduction).
The construction above is a variant of \v{S}ir\v{s}ov \cite{S2}.
Here we choose elimination-on-the-left of all but one variables except the smallest , but we could have done
 elimination-on-the-left of the one greatest variable.
It would have resulted the same Lyndon--Shirshov words but a different evaluation rule.
Altogether, however, what we have here is just noncommutative polynomial elimination with a particularly good accounting.
In general, using Lyndon--Shirshov words to parametrize bases of free Lie algebras may be practical, because
 they lift some burden of organization; they are also the choice of Lothaire \cite{Lot}.
They are advantageous for indexing even if the evaluation rules are very different.
(Cf. Bokut, Chibrikov \cite{BCh}, Chibrikov \cite{Chb}, Walter, Shiri \cite{WS}.)
They fit well to the dimension formula of Witt \cite{W}, cf.~\v{S}ir\v{s}ov \cite{S2}.
(Remark: that reason for changing the ordering compared to the Section \ref{sec:via} was simply
 that we can say `antilexicographical (where tails are smaller)' instead of `antilexicographical but tails are greater'.)
\snewpage
\section{From $\mathrm F_{K}^{\Lie}$ to the basic PBW theorem}\plabel{sec:freefrom}
Here we assume to know that free Lie $K$-algebras are  free $K$-modules in every multigrade separately (MW),
 and that the PBW theorem holds for them.
\begin{prop}
The PBW theorem holds if $\mathfrak g$ is a free $K$-module.
\begin{proof}
Consider a base $\{Z_\lambda:\lambda\in\Lambda\}$ for $\mathfrak g$.
Take $\mathrm F^{\Lie}_K[X_\lambda:\lambda\in\Lambda]$.
Let us extend $\{X_\lambda:\lambda\in\Lambda\}$ by $\{P_\omega:\omega\in\Omega\}$
 obtained from higher multigrades to a basis
\[\mathcal B=\{X_\lambda:\lambda\in\Lambda\}\cup\{P_\omega:\omega\in\Omega\}\]
 of $\mathrm F^{\Lie}_K[X_\lambda:\lambda\in\Lambda]$.
Assume that
\[P_\omega(Z_\lambda:\lambda\in\Lambda)=\sum_{\lambda\in \Lambda}a_{\lambda,\omega}Z_\lambda\]
 ($X_\lambda$ is substituted by $Z_\lambda$ in $P_\omega$). Then let
\[\hat P_\omega=P_\omega-\sum_{\lambda\in \Lambda}a_{\lambda,\omega}X_\lambda.\]
Now
\[\hat{\mathcal B}=\{X_\lambda:\lambda\in\Lambda\}\cup\{\hat P_\omega:\omega\in\Omega\}\]
 is still a basis.
Take any ordering $\leq$ on that; and assume, say, that elements belonging $\Lambda$ precede the ones belonging to $\Omega$.
The elements $\hat P_\omega$  span an ideal $\mathcal I$ in $\mathrm F^{\Lie}_K[X_\lambda:\lambda\in\Lambda]$.
Indeed, they span exactly the kernel of the evaluation map $X_\lambda\mapsto Z_\lambda$.
If $\bo[Z_\lambda,Z_\mu\bo]=\sum_\nu c_{\lambda,\mu}^\nu Z_{\nu}$,
 then $\bo[X_\lambda,X_\mu\bo]-\sum_\nu c_{\lambda,\mu}^\nu X_{\nu}\in \mathcal I$.
Thus   there is a homomorphism
\[\mathfrak g\rightarrow\mathrm  F_K[X_\lambda:\lambda\in\Lambda]/\mathcal I'\]
\[Z_\lambda\mapsto [X_\lambda]_{\mathcal U}/\mathcal I',\]
 where $\mathcal I'$ is the associative ideal generated by the image of $\mathcal I$ through commutator evaluation.
We claim that $\mathcal I'$ is $K$-linearly generated by the elements
\begin{equation}[X_{\lambda_1}\ldots X_{\lambda_n}\hat P_{\omega_1}\ldots\hat P_{\omega_m}]_{\mathcal U} \plabel{eq:gog}
\end{equation}
 $\lambda_1\leq\ldots\leq\lambda_n\leq \omega_1\leq\ldots\leq\omega_m$ such that $n\geq 0$, $m\geq 1$.
Indeed, if we take an arbitrary product of base elements which contains at least one $\hat P_\omega$ and we
 apply the basic rearrangement procedure, then at least one element in any  formal product monomial will be from $\mathcal I$.
A base element from $\{\hat P_\omega:\omega\in\Omega\}$ is either unaffected in a step, or it gets com\-mu\-ta\-ted,
 but then the commutator is a $K$-linear combination of elements of $\{\hat P_\omega:\omega\in\Omega\}$.

Now, the injectivity of $\bo m_\leq$ with respect to $\hat{\mathcal B}$ and
the specific shape of $\mathcal I'$ implies
 that the evaluation map given by
\[\textstyle{\bigotimes_\leq}\mathfrak g\rightarrow\mathrm  F_K[X_\lambda:\lambda\in\Lambda]/\mathcal I'\]
\[Z_{\lambda_1}\otimes\ldots\otimes Z_{\lambda_n}\mapsto  [X_{\lambda_1}\ldots X_{\lambda_n}]_{\mathcal U}/\mathcal I'\]
 ($\lambda_1\leq\ldots\leq\lambda_n$) is injective.
Thus the evaluation map into $\mathcal U\mathfrak g$ must also be injective.
(Remark: Actually, $\mathcal U\mathfrak g\simeq\mathrm   F_K[X_\lambda:\lambda\in\Lambda]/\mathcal I'$
 by universal algebraic reasons.)
\end{proof}
\end{prop}
\snewpage

It seems to be a drawback that we obtained only the basic PBW theorem for free $K$-modules.
This can be remedied as follows. Due to the relatively transparent structure of $\mathrm F^{\Lie}_K[X_\lambda:\lambda\in\Lambda]$,
 one can  define free Lie algebras $\mathrm F^{\Lie}_{K_\Lambda}[X_\lambda:\lambda\in\Lambda]$ with variable coefficient structure.
This means that in multigrade $X_{\lambda_1}^{i_1}\ldots X_{\lambda_s}^{i_s}$ ($i_k\geq1$)
 the coefficient ring is $\botimes^{i_1}(K/I_{\lambda_1})\otimes\ldots\otimes \botimes^{i_s}(K/I_{\lambda_i})$
 ($\simeq K/(I_{\lambda_1}+\ldots +I_{\lambda_s})$).
This has the same monomial structure as $\mathrm F^{\Lie}_{\mathbb Z}[X_\lambda:\lambda\in\Lambda]$; except that
 some multigrades are deselected (where the coefficient ring is $0$), but this makes no essential difference.
This evaluates in the noncommutative polynomial algebra $\mathrm F_{K_\Lambda}[X_\lambda:\lambda\in\Lambda]$, and the
 PBW theorem remains valid, as in every multigrade we have the same evaluation structure as in the free Lie algebra
 with with respect to the appropriate coefficient ring.
But then the arguments of the previous proof can be modified in order to obtain the  basic PBW theorem for sums of cyclic $K$-modules.

\snewpage

\appendix
\section{The Witt--Lazard proof of the global PBW theorems}\plabel{sec:lazard}

Although the classical proofs of the PBW theorem which work for general fields
are quite similar to each other; the approach due to  Witt \cite{W} and
Lazard \cite{L} is characterized (as opposed to Birkhoff \cite{B} and Jacobson \cite{J} or
Cartan, Eilenberg \cite{CE} and Cohn \cite{CC})
by (a) an emphatic appearance of the symmetric group, and (b) a more explicit
description of the ideal structure of the universal factorization.
Thus, it algebraizes the combinatorics quite well.
It allows to formulate the proof of the PBW theorem simultaneously
in (i)  the basic case (sum of cyclic $K$-modules) and (ii) the symmetric case ($\mathbb Q\subset K$).

\textbf{(I) Actions of symmetric groups.}
The symmetric group $\Sigma_n$ acts naturally on $\botimes^n\mathfrak g$
by the presription
\[\sigma*v_1\otimes\ldots\otimes v_n=v_{\sigma^{-1}(1)}\otimes\ldots\otimes v_{\sigma^{-1}(n)}.\]
Then $\botimes^n_{(0)}\mathfrak g$ is generated by $v-\sigma*v$ where $v\in \botimes^n\mathfrak g$, $\sigma\in\Sigma_n$.
Let $W_{k,n}$ denote the permutation $(k\,\,k+1)$ in $\Sigma_n$.
Then it is also true that  $\botimes^n_{(0)}\mathfrak g$ is generated by $v-W_{k,n}*v$ where $v\in \botimes^n\mathfrak g$, $k<n$.
Let us define the $W_{k,n}\bullet:\botimes^n\mathfrak g\rightarrow\botimes^{n-1}\mathfrak g$ by taking a Lie-commutator
between the $k$th and $(k+1)$th positions. So,
\[W_{k,n}\bullet v_1\otimes\ldots\otimes v_n=v_1\otimes\ldots v_{k-1}\otimes[v_k,v_{k+1}]\otimes v_{k+2}\ldots\otimes v_n.\]
We can extend $W_{k,n}* $  and  $W_{k,n}\bullet $ to $\bigotimes\mathfrak g$.
In the first case the action is identity outside $\botimes^n\mathfrak g$, and in the second case it is the zero map.

We define the action $W_{k,n}\diamond  $ as $W_{k,n}* +W_{k,n}\bullet $ (extended sense).
Then $v-W_{k,n}\diamond v $ vanishes if $v\in\botimes^k\mathfrak g$, $n\neq k$.
Let $J^n\mathfrak g$ be the module generated by $v-W_{k,n}\diamond v $, $v\in\botimes^n\mathfrak g$.
We see that $\displaystyle{J\mathfrak g=\sum J^n\mathfrak g}$.
Let us extend $\diamond$ to $\Sigma_n$ as follows.
For $1\in\Sigma_n$  let $1\diamond$ be the identity.
For $\sigma\in\Sigma_n\setminus\{1,W_{1,n},\ldots,W_{n-1,n}\}$ we choose an arbitrary (but fixed)
decomposition $\sigma=W_{a_1,n}\ldots W_{a_s,n}$, and we let $(\sigma\diamond)=(W_{a_1,n}\diamond)\ldots(W_{a_s,n}\diamond)$.
Now,
$\sigma\diamond$ still acts as identity outside $\botimes^n\mathfrak g$, but it does not necessarily define an associative action of $\Sigma_n$.
However, it is not very far from it:
\begin{lemma}\plabel{lem:shallow}
 $W_{k,n}\diamond$ acts trivially on $J^{n-1}\mathfrak g$ (thus invariantly), and
\[(W_{k,n}\diamond )^2=\id\mod J^{n-1}\mathfrak g,\tag{P1}\]
\[(W_{k,n}\diamond )(W_{l,n}\diamond )= (W_{l,n}\diamond )(W_{k,n}\diamond )\mod J^{n-1}\mathfrak g\qquad \text{ if } l-k\geq2,\tag{P2}\]
\[(W_{k,n}\diamond )(W_{k+1,n}\diamond )(W_{k,n}\diamond )= (W_{k+1,n}\diamond )(W_{k,n}\diamond )(W_{k+1,n}\diamond )\mod J^{n-1}\mathfrak g.\tag{P3}\]
\begin{proof}
The triviality property follows from $ J^{n-1}\mathfrak g\subset \botimes^{n-1}\mathfrak g\oplus\botimes^{n-2}\mathfrak g$.
The equalities follow from the identities
\[(W_{k,n}\diamond )(W_{k,n}\diamond )-\id=0\tag{P1'}\]
\begin{multline}
(W_{k,n}\diamond )(W_{l,n}\diamond )- (W_{l,n}\diamond )(W_{k,n}\diamond )=\\
=(\id-(W_{k,n-1}\diamond ))(W_{l,n}\bullet)-(\id-(W_{l-1,n-1}\diamond ))(W_{k,n}\bullet)
\tag{P2'}\end{multline}
if $l-k\geq2$;
\begin{multline}(W_{k,n}\diamond )(W_{k+1,n}\diamond )(W_{k,n}\diamond )-
(W_{k+1,n}\diamond )(W_{k,n}\diamond )(W_{k+1,n}\diamond )=\\
=(\id-(W_{k,n-1}\diamond ))(W_{k,n}\bullet)+(\id-(W_{k,n-1}\diamond ))(W_{k+1,n}\bullet)(W_{k+1,n}*)
 \\  -(\id-(W_{k,n-1}\diamond ))(W_{k+1,n}\bullet);\tag{P3'}\end{multline}
which, checked against $v_1\otimes\ldots\otimes v_n$, follow from  the Lie-identities.
\end{proof}
\end{lemma}

\begin{cor}
$\diamond$ extends to an associative action of $\Sigma_n$ modulo $J^{n-1}\mathfrak g$.
\begin{proof}
In (P1)--(P3) we recognize the semigroup presentation of $\Sigma_n$ based on $W_{k,n}$ (Cf.~ Dickson \cite{DD}, P. 2, Ch. XIII).
The relations are satisfied according to the previous lemma, thus the action descends to $\Sigma_n$.
\end{proof}
\end{cor}

\textbf{(II) The tensorial splittings.} We define  the forgetting map
 $\mathrm f^n:\Sigma_n\otimes   \botimes^n\mathfrak g \rightarrow \botimes^n\mathfrak g$
 such that
\[\mathrm f^n(\sigma\otimes v_1\otimes\ldots\otimes v_n)=v_1\otimes\ldots\otimes v_n;\]
 and we define the evaluation map
 $\mathrm e^n:\Sigma_n\otimes  \botimes^n\mathfrak g \rightarrow \botimes^n\mathfrak g$
 such that
\[\mathrm e^n(\sigma\otimes v_1\otimes\ldots\otimes v_n)=\sigma*v_1\otimes\ldots\otimes v_n.\]

In case (i), we take a basis a $\{g_\alpha \,:\,\alpha\in A\}$, and introduce an ordering $\leq$ on $A$.
We define  $ \eta_n:\botimes^n\mathfrak g\rightarrow\Sigma_n\otimes  \botimes^n\mathfrak g$ such that
\[\eta_n(g_{\alpha_1}\otimes\ldots \otimes g_{\alpha_n})=\sigma\otimes g_{\alpha_1}\otimes\ldots\otimes g_{\alpha_n},\]
 where $\alpha_{\sigma^{-1}(1)}\leq\ldots\leq\alpha_{\sigma^{-1}(n)}$ and $i< j$, $\alpha_i=\alpha_j$ implies
 $\sigma(i)<\sigma(j)$.
I.~e.~$\sigma$ is the permutation which orders $\alpha_1,\ldots,\alpha_n$ with the least number of involutions.

In case (ii), we define
\[\eta_n(g_{1}\otimes\ldots\otimes g_{n})=\frac1{n!}\sum_{\sigma\in\Sigma_n}\sigma\otimes g_{1}\otimes\ldots \otimes g_{n}.\]
Then
\[\mathrm f^{n}\circ\eta_n=\id.\]
It easy to see from the definitions that
\begin{equation}\mathrm e^{n}\circ\eta_n(v)=\mathrm e^{n}\circ\eta_n(\sigma*v)\plabel{eq:invar}\end{equation}
 for any $\sigma\in\Sigma_n$, $v\in\botimes^{n}\mathfrak g$.
This is the same thing to say  as $\mathrm e^{n}\circ\eta_n\circ(\mathrm f^n-\mathrm e^n)=0$.
Then
\[\mathrm e^{n}\circ\eta_n\text{ is an idempotent}.\]
Indeed, 
 $\mathrm e^{n}\circ\eta_n\circ(\id-\mathrm e^{n}\circ\eta_n)=\mathrm e^{n}\circ\eta_n\circ(\mathrm f^n-\mathrm e^{n})\circ\eta_n=0$.

This idempotence yields the direct sum decomposition
\[\botimes^n\mathfrak g =\underbrace{\botimes^n_{\eta}\mathfrak g}_{\im \mathrm e^{n}\circ\eta_n }\oplus
\underbrace{\qquad\botimes^n_{(0)}\mathfrak g\qquad}_{\im \id -\mathrm e^{n}\circ\eta_n =\ker \mathrm e^{n}\circ\eta_n },\]
where the first factor is named so by definition, and regarding the identification of the second factor we note that
$\im (\id-\mathrm e^{n}\circ\eta_n)=\im \mathrm (f^n-\mathrm e^{n})\circ\eta_n\subset\botimes^n_{(0)}\mathfrak g\equiv\im \mathrm f^n-\mathrm e^{n}\subset \ker \mathrm e^{n}\circ\eta_n$.
(Note that very little happens here.)

\textbf{(III) The PBW splittings.}
Note that in case (i), $\botimes_{\eta}\mathfrak g=\botimes_\leq\mathfrak g$; and in case (ii), $\botimes_{\eta}\mathfrak g=\botimes_\Sigma \mathfrak g$.
Thus, the statement of the PBW theorem is that $\bigotimes_\eta\mathfrak g$ and $J\mathfrak g$
do not intersect each other (and, in fact, they are complementer spaces in $\bigotimes\mathfrak g$).


We define the evaluation map
$\mathrm e\mathrm e^{n}:\Sigma_n\otimes  \botimes^n\mathfrak g
\rightarrow \botimes^n\mathfrak g\oplus  \botimes^{n-1}\mathfrak g$
such that
\[\mathrm e\mathrm e^{n}(\sigma\otimes v_1\otimes\ldots\otimes v_n)=\sigma\diamond v_1\otimes\ldots\otimes v_n.\]
\snewpage
\begin{lemma}\plabel{lem:lazardkey} For $v\in\botimes^n\mathfrak g$,
\begin{equation}
v-W_{k,n}\diamond v=(\id- \mathrm e\mathrm e^{n}\circ\eta_n)(v-W_{k,n}*v) \mod J^{n-1}\mathfrak g.
\plabel{eq:lazardkey}
\end{equation}
\begin{proof}
Let us note that $\sigma\diamond$, $\sigma\in\Sigma_n$, acts trivially on $W_{k,n}\diamond v-W_{k,n}* v  \in\botimes^{n-1}\mathfrak g$; so
$(\id-\sigma\diamond)(W_{k,n}\diamond v-W_{k,n}* v)=0$. This implies   $(\id-\sigma\diamond)(W_{k,n}* v)=(\id-\sigma\diamond)(W_{k,n}\diamond v)$.

Case (i): Assume that $v=g_{\alpha_1}\otimes\ldots \otimes g_{\alpha_n}$.
If $\alpha_k=\alpha_{k+1}$, then $W_{k,n}*v=v$, $W_{k,n}\diamond v=v$,
\[(\id- \mathrm e\mathrm e^{n}\circ\eta_n)(v-W_{k,n}* v)=0=v-W_{k,n}\diamond v.\]
If $\alpha_k\neq\alpha_{k+1}$,  $\eta_n(v)=\sigma \otimes v$, then $\eta_n(W_{k,n}*v)=(\sigma W_{k,n} ) \otimes (W_{k,n}*v)$. Thus
\begin{align}\notag
(\id- \mathrm e\mathrm e^{n}\circ\eta_n)(v-W_{k,n}* v)&=(\id- \mathrm e\mathrm e^{n}\circ\eta_n)(v)-(\id- \mathrm e\mathrm e^{n}\circ\eta_n)(W_{k,n}* v)
\\\notag&=
(\id-\sigma\diamond)v-(\id-\sigma W_{k,n}\diamond)(W_{k,n}*v)
\\\notag&=(\id-\sigma\diamond)v-(\id-\sigma W_{k,n}\diamond)(W_{k,n}\diamond v)
\\\notag&\hspace{-7.7mm}
\stackrel{\mod J^{n-1}\mathfrak g }{=}(v-\sigma\diamond v)-( W_{k,n}\diamond v- \sigma\diamond v)
\\\notag&=v-W_{k,n}\diamond v.
\end{align}

Case (ii):
\begin{align}\notag
(\id- \mathrm e\mathrm e^{n}\circ\eta_n)(v-W_{k,n}* v)&=
(\id- \mathrm e\mathrm e^{n}\circ\eta_n)(v)-(\id- \mathrm e\mathrm e^{n}\circ\eta_n)(W_{k,n}* v)
\\\notag&=
\frac1{n!}\sum_{\sigma\in\Sigma_n}(\id-\sigma\diamond)v-\frac1{n!}\sum_{\sigma\in\Sigma_n}(\id-\sigma\diamond)(W_{k,n}*v)
\\\notag&=
\frac1{n!}\sum_{\sigma\in\Sigma_n}(\id-\sigma\diamond)v-\frac1{n!}\sum_{\sigma\in\Sigma_n}(\id-\sigma\diamond)(W_{k,n}\diamond)v
\\\notag&\hspace{-7.6mm}\stackrel{\mod J^{n-1}\mathfrak g}{=}
v-\frac1{n!}\sum_{\sigma\in\Sigma_n}\sigma\diamond v-W_{k,n}\diamond v+\frac1{n!}\sum_{\sigma\in\Sigma_n} \sigma W_{k,n}\diamond v
\\\notag&=v-W_{k,n}\diamond v.\qedhere
\end{align}
\end{proof}
\end{lemma}
\begin{commentx}
\begin{remark}\plabel{rem:lazardkey}
The elements  $v-\sigma\diamond v$  with $v\in \botimes^n\mathfrak g$, $\sigma\in\Sigma_n$ still
generate only $J^{n}\mathfrak g$.
Indeed, if the canonical decomposition is $\sigma=W_{a_1,n}\ldots W_{a_s,n}$, then
\begin{align}\notag
v-\sigma\diamond v =&\sum_{i=1}^s  (\id- (W_{a_i,n}\diamond)) (W_{a_{i+1},n}\diamond) \ldots (W_{a_s,n}\diamond)v\\
=&\sum_{i=1}^s  (\id- (W_{a_i,n}\diamond)) (W_{a_{i+1},n}*) \ldots (W_{a_s,n}*)v\in J^{n}\mathfrak g.
\notag\end{align}

Using similar arguments, $W_{k,n}$ can be replaced by an arbitrary $\sigma\in\Sigma_n$ in equation \eqref{eq:lazardkey}.
Actually, the discussion yields constructive maps $\mathrm h_n:\Sigma_n\otimes\botimes^{n}\mathfrak g\rightarrow W\otimes\botimes^{n-1}\mathfrak g$
such that
\[v-\sigma\diamond v=(\id- \mathrm e\mathrm e^{n}\circ\eta_n)(v-\sigma*v)+(\mathrm f^{n-1}- \mathrm e^{n-1})\circ \mathrm h_n(\sigma,v) .\]

If we consider only the $\botimes^n\mathfrak g$-part, then this yields
 $v-\sigma* v=(\id- \mathrm e^{n}\circ\eta_n)(v-\sigma*v)$, which simplifies to $\mathrm e^{n}\circ\eta_n(v)=\mathrm e^{n}\circ\eta_n(\sigma*v)$, cf.~equation \eqref{eq:invar}.

(This remark is not needed to the proof.)\qedremark
\end{remark}
\end{commentx}
Let $J_{\eta}^n\mathfrak g$ be the image of $\botimes^n_{(0)}\mathfrak g$ under   $\id- \mathrm e\mathrm e^{n}\circ\eta_n=(\mathrm f^n-\mathrm e\mathrm e^n)\circ\eta_n$.
(Strictly speaking, $J_{\eta}^n\mathfrak g$  depends not only on $\eta$ but also on $\diamond$.)
Note that any element $w\in J_{\eta}^n\mathfrak g$ can be reconstructed from its projection to $\botimes^n\mathfrak g$ which is in $\botimes^n_{(0)}\mathfrak g$.
Indeed, the projection of $v-W_{k,n}\diamond v$ is $v-W_{k,n}* v$, and the projection of $J^{n-1}\mathfrak g$ is $0$; and formula \eqref{eq:lazardkey} implies
$w= (\id- \mathrm e\mathrm e^{n}\circ\eta_n)(\proj_{\,\bigotimes^n\mathfrak g}w)$.
\begin{commentx}
\begin{remark}\plabel{rem:split}
One can show that $J_{\eta}^n\mathfrak g\subset J^n\mathfrak g$.
Indeed, the LHS is  the image of $\botimes^n_{(0)}\mathfrak g$ under   $(\mathrm f^n-\mathrm e\mathrm e^n)\circ\eta_n$,
while the RHS is the image of $\mathrm f^n-\mathrm e\mathrm e^n$ (cf.~the beginning of the previous Remark).
Now, $J^n\mathfrak g$ projects to $\botimes^n_{(0)}\mathfrak g$,
thus $(\id- \mathrm e\mathrm e^{n}\circ\eta_n)\circ\proj_{\,\bigotimes^n\mathfrak g}$ yields an idempotent on
$J^n\mathfrak g$.
It is straightforward  to see from Lemma \ref{lem:lazardkey} that the corresponding inner direct sum decomposition is
\[J^n\mathfrak g=\underbrace{\qquad J^n_\eta\mathfrak g \qquad}_{\im (\id- \mathrm e\mathrm e^{n}\circ\eta_n)\circ\proj_{\,\bigotimes^n\mathfrak g}}
\oplus \underbrace{(J^n\mathfrak g\cap  J^{n-1}\mathfrak g)}_{\ker (\id- \mathrm e\mathrm e^{n}\circ\eta_n)\circ\proj_{\,\bigotimes^n\mathfrak g}}.\]
(E. g. $[x,y]\otimes z+[y,z]\otimes x+[z,x]\otimes y-z\otimes[x,y]-x\otimes[y,z]-y\otimes[z,x]\in J^3\mathfrak g\cap J^2\mathfrak g$.)

(This remark is not needed to the proof.)
\qedremark
\end{remark}
\end{commentx}

Let $J^{\leq n}\mathfrak g=J^{0}\mathfrak g+\ldots+J^{n}\mathfrak g$.

\begin{cor}\plabel{cor:split}
The following inner direct sum decompositions hold:

(a) $J^{\leq n}\mathfrak g=J^{\leq n-1}\mathfrak g\oplus  J_{\eta}^n\mathfrak g$;

(b) $J^{\leq n}\mathfrak g= J_{\eta}^0\mathfrak g\oplus\ldots\oplus J^{n}_{\eta}\mathfrak g$;

(c) $J\mathfrak g= J_{\eta}^0\mathfrak g\oplus\ldots\oplus J^{n}_{\eta}\mathfrak g\oplus\ldots$;

(d) $\botimes^{\leq n}\mathfrak g=\botimes^{\leq n-1}\mathfrak g\oplus \botimes^{n}\mathfrak g=\botimes^{\leq n-1}\mathfrak g\oplus \botimes^{n}_\eta\mathfrak g\oplus J^n_\eta\mathfrak g$;

(e)
$\underbrace{\botimes^0\mathfrak g \oplus\ldots\oplus\botimes^n\mathfrak g}_{\textstyle{\bigotimes}^{\leq n}\mathfrak g} =
\underbrace{\botimes^0_\eta\mathfrak g \oplus\ldots\oplus\botimes^n_\eta\mathfrak g}_{\textstyle{\bigotimes}^{\leq n}_\eta\mathfrak g} \oplus
\underbrace{   J_{\eta}^0\mathfrak g\oplus\ldots\oplus J^{n}_{\eta}\mathfrak g}_{\textstyle{J^{\leq n}\mathfrak g}}$;

(f) $\textstyle{\bigotimes}\mathfrak g =\textstyle{\bigotimes}_\eta\mathfrak g\oplus J\mathfrak g$.
\\
Remark: Here $\botimes^0\mathfrak g=\botimes^0_\eta\mathfrak g=K$, $J^0\mathfrak g=J^0_\eta\mathfrak g=0$, $J^1\mathfrak g=J^1_\eta\mathfrak g=0$.
\begin{proof}
(a) Lemma \ref{lem:lazardkey} implies $J^{n-1}\mathfrak g+J^{n}\mathfrak g=J^{n-1}\mathfrak g+J_\eta^{n}\mathfrak g$.
Adding $J^{\leq n-2}\mathfrak g$ to both sides, we obtain $J^{\leq n}\mathfrak g =J^{\leq n-1}\mathfrak g+J_\eta^{n}\mathfrak g$.
The two factors in the sum must be disjoint, as projected to $\botimes^n\mathfrak g$,  the first factor projects to $0$, while the
second factor projects faithfully.

(b) follows from (a) inductively.

(c) follows from (b) by taking increasing unions.

(d) The first equality in obvious. The second one follows from the fact that on the RHS,
 as projected to $\botimes^n\mathfrak g$, the second factor projects to $\botimes^{n}_\eta\mathfrak g$ faithfully, and
  $J^n_\eta\mathfrak g$ projects to
$\botimes_{(0)}^n\mathfrak g$ faithfully.

(e) follows from (d) inductively, the labeling uses (b).

(f) follows from (e) by taking increasing unions.
\end{proof}
\end{cor}
In particular, we find the statement of the PBW theorem in (f).

Specified to the basic case (i), the proof given above is certainly over-algebraized compared to the one of Jacobson \cite{J}
 (but which is very similar in content).
However, it conveys the message that the point is not the cleverness of the reduction process
 but the control over the factorization ideal, which is permitted by the shallowness
 of the ambiguities, embodied in Lemma \ref{lem:shallow}.

The particular statement  ${\textstyle{\bigotimes}^{\leq n}\mathfrak g} =
 {\textstyle{\bigotimes}^{\leq n}_\eta\mathfrak g} \oplus
 {\textstyle{J^{\leq n}\mathfrak g}}$ from (e) appears explicitly in Lazard \cite{L}
 (and it is also implicit in Jacobson \cite{J}).
This fact combined with the observation that the possible collapse of the $\bo m^{(n)}$ is a finitistic matter,
 is the starting point toward the case when $K$ is a principal ideal domain or a Dedekind domain,
 see Lazard \cite{L}, Cartier \cite{C}.

There are similar proofs which concentrate not on the factorization structure but on constructing
 the regular representation of the enveloping algebra.
Such one is the proof of Cartan, Eilenberg \cite{CE} in the basic case (i).
In the symmetric case (ii), this is done by Cohn \cite{CC};
 another (terse) proof is given in Deligne, Morgan, Bernstein \cite{DMB}.

Quite different proof is the one of Milnor, Moore \cite{MiMo} (over fields)
 which uses the coalgebra structure of the universal enveloping algebra.
(We also make a similar argument but leaning on the Magnus--Witt theorems.)

The combinatorial-in-ring content is conceptualized in
Bergman \cite{Ber}.
The general approach of this paper is also somewhat combinatorial but it takes place ``outside of
 the enveloping algebra of the original Lie algebra'' in the enveloping algebra of
 the free Lie algebras.

Graded versions of the PBW theorem as such provide no particular challenges.
In general, however, theorems of PBW type can be drawn in several contexts, with various
difficulty,
cf. Shepler, Witherspoon \cite{SW}.  

\snewpage

\end{document}